\documentstyle[12pt,amscd]{amsart}
\newtheorem{thm}[equation]{Theorem}
\numberwithin{equation}{section}
\newtheorem{cor}[equation]{Corollary}

\newtheorem{rmk}[equation]{Remark}

\newtheorem{lem}[equation]{Lemma}

\newtheorem{diag}[equation]{Diagram}
\newtheorem{prop}[equation]{Proposition}

\begin{document}
\raggedbottom \voffset=-.7truein \hoffset=0truein \vsize=8truein
\hsize=6truein \textheight=8truein \textwidth=6truein
\baselineskip=18truept
\def\mapright#1{\ \smash{\mathop{\longrightarrow}\limits^{#1}}\ }
\def\mapleft#1{\smash{\mathop{\longleftarrow}\limits^{#1}}}
\def\mapup#1{\Big\uparrow\rlap{$\vcenter {\hbox {$#1$}}$}}
\def\mapdown#1{\Big\downarrow\rlap{$\vcenter {\hbox {$\ssize{#1}$}}$}}
\def\mapne#1{\nearrow\rlap{$\vcenter {\hbox {$#1$}}$}}
\def\mapse#1{\searrow\rlap{$\vcenter {\hbox {$\ssize{#1}$}}$}}
\def\mapr#1{\smash{\mathop{\rightarrow}\limits^{#1}}}
\def\ss{\smallskip}
\def\vp{v_1^{-1}\pi}
\def\at{{\widetilde\alpha}}
\def\sm{\wedge}
\def\la{\langle}
\def\ra{\rangle}
\def\on{\operatorname}
\def\spin{\on{Spin}}
\def\kbar{{\overline k}}
\def\qed{\quad\rule{8pt}{8pt}\bigskip}
\def\ssize{\scriptstyle}
\def\a{\alpha}
\def\bz{{\Bbb Z}}
\def\im{\on{im}}
\def\ct{\widetilde{C}}
\def\ext{\on{Ext}}
\def\sq{\on{Sq}}
\def\eps{\epsilon}
\def\ar#1{\stackrel {#1}{\rightarrow}}
\def\br{{\bold R}}
\def\bC{{\bold C}}
\def\bA{{\bold A}}
\def\bB{{\bold B}}
\def\bD{{\bold D}}
\def\bh{{\bold H}}
\def\bQ{{\bold Q}}
\def\bP{{\bold P}}
\def\bx{{\bold x}}
\def\bo{{\bold{bo}}}
\def\si{\sigma}
\def\Ebar{{\overline E}}
\def\dbar{{\overline d}}
\def\Sum{\sum}
\def\tfrac{\textstyle\frac}
\def\tb{\textstyle\binom}
\def\Si{\Sigma}
\def\w{\wedge}
\def\equ{\begin{equation}}
\def\b{\beta}
\def\G{\Gamma}
\def\g{\gamma}
\def\k{\kappa}
\def\psit{\widetilde{\Psi}}
\def\tht{\widetilde{\Theta}}
\def\psiu{{\underline{\Psi}}}
\def\thu{{\underline{\Theta}}}
\def\aee{A_{\text{ee}}}
\def\aeo{A_{\text{eo}}}
\def\aoo{A_{\text{oo}}}
\def\aoe{A_{\text{oe}}}
\def\fbar{{\overline f}}
\def\endeq{\end{equation}}
\def\sn{S^{2n+1}}
\def\zp{\bold Z_p}
\def\A{{\cal A}}
\def\P{{\cal P}}
\def\cj{{\cal J}}
\def\zt{{\bold Z}_2}
\def\bs{{\bold s}}
\def\bof{{\bold f}}
\def\bq{{\bold Q}}
\def\bu{{\bold u}}
\def\be{{\bold e}}
\def\Hom{\on{Hom}}
\def\ker{\on{ker}}
\def\coker{\on{coker}}
\def\da{\downarrow}
\def\colim{\operatornamewithlimits{colim}}
\def\zphat{\bz_2^\wedge}
\def\io{\iota}
\def\Om{\Omega}
\def\u{{\cal U}}
\def\e{{\cal E}}
\def\exp{\on{exp}}
\def\wbar{{\overline w}}
\def\xbar{{\overline x}}
\def\ybar{{\overline y}}
\def\zbar{{\overline z}}
\def\ebar{{\overline e}}
\def\nbar{{\overline n}}
\def\rbar{{\overline r}}
\def\et{{\widetilde E}}
\def\ni{\noindent}
\def\coef{\on{coef}}
\def\den{\on{den}}
\def\lcm{\on{l.c.m.}}
\def\vi{v_1^{-1}}
\def\ot{\otimes}
\def\psibar{{\overline\psi}}
\def\mhat{{\hat m}}
\def\exc{\on{exc}}
\def\ms{\medskip}
\def\ehat{{\hat e}}
\def\etao{{\eta_{\text{od}}}}
\def\etae{{\eta_{\text{ev}}}}
\def\dirlim{\operatornamewithlimits{dirlim}}
\def\gt{\widetilde{L}}
\def\lt{\widetilde{\lambda}}
\def\st{\widetilde{s}}
\def\ft{\widetilde{f}}
\def\sgd{\on{sgd}}
\def\lfl{\lfloor}
\def\rfl{\rfloor}
\def\ord{\on{ord}}
\def\gd{{\on{gd}}}
\def\rk{{{\on{rk}}_2}}
\def\nbar{{\overline{n}}}
\def\lg{{\on{lg}}}
\def\N{{\Bbb N}}
\def\Z{{\Bbb Z}}
\def\Q{{\Bbb Q}}
\def\R{{\Bbb R}}
\def\C{{\Bbb C}}
\def\l{\left}
\def\r{\right}
\def\mo{\on{mod}}
\def\tmf{\on{tmf}}
\def\vexp{v_1^{-1}\exp}
\def\notimm{\not\subseteq}
\def\Remark{\noindent{\it  Remark}}
\title[Nonimmersions implied by tmf, revisited]
{Nonimmersions of $RP^n$ implied by $\tmf$, revisited}
\author{Donald M. Davis}
\address{Lehigh University\\Bethlehem, PA 18015, USA}
\email{dmd1@@lehigh.edu}

\author{Mark Mahowald}
\address{Northwestern University\\Evanston, IL 60208, USA}
\email{mark@@math.northwestern.edu}
\date{April 5, 2007}

\keywords{immersion, projective space, elliptic cohomology}
\thanks{We thank Steve Wilson for causing us to take a look
at these matters.}
\subjclass[2000]{57R42, 55N20.}

\begin{abstract} In a 2002 paper, the authors and Bruner used
the new spectrum $\tmf$ to obtain some new nonimmersions of real
projective spaces. In this note, we complete/correct two oversights
in that paper. 

The first is to note that in that paper a general
nonimmersion result was stated which yielded new nonimmersions for $RP^n$ with
$n$ as small
as 48, and yet it was stated there that the first new result occurred
when $n=1536$. Here we give a simple proof of those overlooked results.

Secondly, we fill in a gap in the proof of the 2002 paper. There it was 
claimed that an axial map $f$ must satisfy
$f^*(X)=X_1+X_2$.  We realized recently that this is not clear.
 However, here we show that it is true up multiplication by a unit
 in the appropriate ring, and so we retrieve all the nonimmersion results
 claimed in \cite{BDM}.
 
  Finally, we present a complete determination  
 of $\tmf^{8*}(RP^\infty\times RP^\infty)$ and $\tmf^*(CP^\infty\times CP^\infty)$ in positive
 dimensions. 
\end{abstract}
\maketitle
\section{Introduction}\label{intro}
In \cite{BDM}, the authors and Bruner described a proof of the following theorem, along with some
additional nonimmersion results.
\begin{thm}\label{BDMthm}$(\cite[1.1]{BDM})$ Assume that $M$ is divisible by the smallest
$2$-power greater than or equal to $h$.
\begin{itemize}
\item If $\a(M)=4h-1$, then $P^{8M+8h+2}$ cannot be immersed in $(\not\subseteq)$
$\R^{16M-8h+10}$.
\item If $\a(M)=4h-2$, then $P^{8M+8h}\not\subseteq\R^{16M-8h+12}$.
\end{itemize}
\end{thm}
\ni Here and throughout, $\a(M)$ denotes the number of 1's in the binary expansion of $M$, and $P^n$
denotes real projective space.

In \cite{BDM}, the theorem is followed by a comment that this is new provided $\a(M)\ge 6$, i.e., $h\ge2$, and the
first new result occurs for $P^{1536}$. In this note, we point out that \ref{BDMthm} is valid when
$h=1$, and these results are new when $M$ is even, including new nonimmersions of $P^n$ for $n$ as small as
56. A remark in \cite[p.66]{BDM} that the nonimmersions when $h=1$ were implied by earlier work of the authors
was incorrect.
Letting $h=1$ in \ref{BDMthm}, we have the following result.
\begin{cor}\label{h=1} a. If $\a(M)=3$, then $P^{8M+10}\notimm\R^{16M+2}$.

b. If $\a(M)=2$, then $P^{8M+8}\notimm\R^{16M+4}$.
\end{cor}

Part (a) is new when $M$ is even. It is 2 better than the previous best result, proved in \cite{AD}, and the nonembedding result that it implies is also new, 1 better than the previous best, proved in \cite{As}.
In \cite{immtable}, a table of known nonimmersions, immersions, nonembeddings, and embeddings of $P^n$
is presented, arranged according to $n=2^i+d$ with $0\le d<2^i$ and $d<64$. Part (a) enters the table with a new result
for $d=58$, applying first to $P^{122}$.

If $M$ is even, \ref{h=1}.b is new, 1 better than the previous best result, of \cite{DZ}, and the nonembedding
result implied is also new. It enters \cite{immtable} at $d=24$ and $40$, with a new result for $P^n$ with $n$ as small as 56. The result of \ref{h=1}.b with $M=2^i+1$ was also proved very recently by Kitchloo and Wilson in \cite{KW}. This result for $P^{2^k+16}$,
2 better than the previous result of \cite{AD} and also new as a nonembedding, enters \cite{immtable} at $d=16$, and applies for $n$ as small as 48.

In Section \ref{h=1sec}, we present a self-contained proof of Corollary \ref{h=1}. The primary reason for doing this,
which amounts to a reproof of part of \cite[1.1]{BDM}, is 
that the proof of the general case in \cite{BDM} requires some extremely elaborate arguments
and calculations. Our proof here, which is just for the case $h=1$, is much more comprehensible.

The proof in \cite{BDM} contained an oversight which we shall correct here. The argument there was that
an immersion of $RP^n$ in $\R^{n+k}$ implies existence of an axial map $P^n\times P^m@>f>> P^{m+k}$
for an appropriate value of $m$, and obtains a contradiction for certain $n$, $m$, and $k$ by consideration of $\tmf^*(f)$.  Here $\tmf$ is the spectrum of topological modular forms, which was discussed in \cite{BDM}. A class $X\in \tmf^8(P^n)$ was described, along with $X_1=X\times 1$ and $X_2=1\times X$ in
$\tmf^8(P^n\times P^m)$. It was asserted that $f^*(X)=X_1+X_2$, and a contradiction obtained by showing
that, for certain values of the parameters, we might have $X^\ell=0$ but $(X_1+X_2)^\ell\ne0$.
We recently realized that it is conceivable that $f^*(X)$ might contain other terms coming from $\tmf^8(P^n\w P^m)$. 

In Section \ref{tPPsec} (see Theorem \ref{tPPthm}) we perform a complete
calculation of $\tmf^*(P^\infty\times P^\infty)$ in positive gradings divisible by 8, 
and in Section \ref{axialsec} we use it to show that effectively $f^*(X)=\bu(X_1+X_2)$, where $\bu$
is a unit in $\tmf^*(P^\infty\times P^\infty)$, which enables us to retrieve
all the nonimmersions of \cite{BDM}.

In Section \ref{CPsec}, we compute $\tmf^*(CP^\infty\times CP^\infty)$ in positive gradings.
The original purpose of doing this was, prior to our obtaining the argument of Section \ref{axialsec},
to see whether we might mimic the argument of \cite{As1} and \cite{Ann} to conclude that if $f$ is an
axial map, then $f^*(X)$ might necessarily equal $u(X_1-X_2)$, where $u$ is a unit in
$\tmf^*(CP\times CP)$. This approach to retrieving the nonimmersions of \cite{BDM} did not yield
the desired result, but the later approach given in Section \ref{axialsec} did. Nevertheless
the nice result for $\tmf^*(CP^\infty\times CP^\infty)$
obtained in Theorem \ref{CPxCP} should be of independent interest.

\section{Proof of Corollary \ref{h=1}}\label{h=1sec}
We begin by proving \ref{h=1}.a. The following standard reduction goes back at least to \cite{Ja}.
If $P^{8M+10}\subseteq\R^{16M+2}$, then gd$((2^{L+3}-8M-11)\xi_{8M+10})\le 8M-8$, hence this bundle has $(2^{L+3}-16M-3)$ linearly independent sections, and thus there is an axial map
$$P^{8M+10}\times P^{2^{L+3}-16M-4}@>f>> P^{2^{L+3}-8M-12}.$$
The bundle here is the stable normal bundle, $L$ is a sufficiently large integer, and gd refers to geometric dimension.
Let $X$, $X_1$, and $X_2$ be elements of $\tmf^8(-)$ described in \cite{BDM} and also in Section \ref{intro}. In Section \ref{axialsec}, we will show that we may assume that $f^*(X)=X_1+X_2$, as was done in \cite{BDM}, since this is true
up to multiplication by a unit. Since
$\tmf^{2^{L+3}-8M-8}(P^{2^{L+3}-8M-12})=0$, we have
$$0=f^*(0)=f^*(X^{2^L-M-1})=(X_1+X_2)^{2^L-M-1}\in \tmf^{2^{L+3}-8M-8}(P^{8M+10}\times P^{2^{L+3}-16M-4}).$$
Expanding, we obtain $\binom{2^L-M-1}{M+1}X_1^{M+1}X_2^{2^L-2M-2}+\binom{2^L-M-1}MX_1^MX_2^{2^L-2M-1}$
as the only terms which are possibly nonzero.
Next we note that, with all $u$'s representing odd integers,
$$\tbinom{2^L-M-1}{M+1}=u_1\tbinom{2M+1}{M+1}=2^{\a(M)-\nu(M+1)}u_2=2^{3-\nu(M+1)}u_2,$$ where we have used $\a(M)=3$ at the last step. Here and throughout, $\nu(2^eu)=e$. Similarly, $\binom{2^L-M-1}M=u_3\binom{2M}M=2^{\a(M)}u_4=2^3u_4$.
Thus an immersion implies that in $\tmf^{2^{L+3}-8M-8}(P^{8M+10}\times P^{2^{L+3}-16M-4})$, we have
\begin{equation}\label{class}2^{3-\nu(M+1)}u_2X_1^{M+1}X_2^{2^L-2M-2}+2^3u_4X_1^MX_2^{2^L-2M-1}=0.\end{equation}

We recall \cite[2.6]{BDM}, which states that there is an equivalence of spectra $P_{b+8}^{k+8}\w \tmf\simeq \Sigma^8P_b^k\w \tmf $.
Combining this with duality, we obtain $\tmf^{8M+8}(P^{8M+10})\approx \tmf_{-1}(P_{-3})\approx\Z/8$, and so $8X_1^{M+1}X_2^{2^L-2M-2}=0$. Here and throughout, $P_n=P_n^\infty=RP^\infty/RP^{n-1}$.
Similarly $\tmf^{2^{L+3}-16M-8}(P^{2^{L+3}-16M-4})\approx \tmf_{7}(P_{3})\approx\Z/16$, and hence
$16X_1^MX_2^{2^L-2M-1}=0$. Duality also implies
$$\tmf^{2^{L+3}-8M-8}(P^{8M+10}\times P^{2^{L+3}-16M-4})\approx \tmf_{14}(P_{-3}\w P_3).$$
Calculations such as $E_2(\tmf_*(P_{-3}\w P_3))$, the $E_2$-term of the Adams spectral sequence (ASS), were made by Bruner's minimal-resolution computer programs in our work on
\cite{BDM}. This one is in a small enough range to actually do by hand. The result is given in Diagram
\ref{diag1}.

\begin{center}
\begin{minipage}{6.5in}
\begin{diag}\label{diag1}{$E_2(\tmf_*(P_{-3}\w P_3))$, $*\le15$}
\begin{center}
\begin{picture}(440,180)
\def\mp{\multiput}
\def\elt{\circle*{3}}
\put(75,20){\line(1,0){310}}
\put(78,5){$0$}
\put(138,5){$3$}
\put(218,5){$7$}
\put(298,5){$11$}
\put(376,5){$15$}
\mp(80,20)(20,0){5}{\elt}
\mp(100,20)(140,0){2}{\line(1,1){40}}
\mp(120,40)(20,20){2}{\elt}
\mp(140,40)(20,0){3}{\elt}
\mp(140,20)(140,0){2}{\line(0,1){40}}
\mp(160,20)(80,40){2}{\line(1,1){20}}
\mp(240,60)(20,20){2}{\elt}
\mp(200,20)(0,20){2}{\elt}
\put(200,20){\line(0,1){20}}
\put(202,38){\elt}
\mp(220,20)(0,20){4}{\elt}
\mp(220,20)(140,0){2}{\line(0,1){60}}
\mp(360,20)(0,20){4}{\elt}
\mp(240,40)(20,20){2}{\elt}
\mp(220,20)(40,80){2}{\line(1,1){40}}
\mp(260,100)(20,20){3}{\elt}
\mp(300,100)(0,20){3}{\elt}
\put(300,100){\line(0,1){40}}
\mp(240,20)(20,20){3}{\elt}
\mp(280,20)(0,20){2}{\elt}
\put(283,20){\elt}
\mp(300,60)(20,0){2}{\elt}
\put(320,40){\line(1,1){40}}
\mp(320,40)(20,20){2}{\elt}
\put(357,20){\line(0,1){38}}
\mp(357,20)(0,19){3}{\elt}
\put(319,20){\line(1,1){38}}
\mp(319,20)(19,19){2}{\elt}
\mp(363,80)(0,20){3}{\elt}
\put(363,80){\line(0,1){40}}
\put(363,80){\line(1,1){17}}
\mp(380,97)(0,-19){3}{\elt}
\put(380,97){\line(0,-1){38}}
\put(383,20){\line(0,1){140}}
\mp(383,20)(0,20){8}{\elt}
\end{picture}
\end{center}
\end{diag}
\end{minipage}
\end{center}

The $\Z/8\oplus\Z/16$ arising from filtration 0 in grading 14 in \ref{diag1} is not hit by a differential
from the class in $(15,0)$ because, as explained in the last paragraph of page 54 of \cite{BDM}, the class
in $(15,0)$ corresponds to an easily-constructed nontrivial map. The monomials $X_1^{M+1}X_2^{2^L-2M-2}$ and
$X_1^MX_2^{2^L-2M-1}$ are detected in mod-2 cohomology, and so their duals emanate from filtration 0.
We saw in the previous paragraph that 8 and 16, respectively, annihilate these monomials, and hence also
their duals. Since the chart shows that the subgroup of $\tmf_{14}(P_{-3}\w P_3)$ generated by classes of filtration
0 is $\Z/8\oplus\Z/16$, we conclude that 8 and 16, respectively, are the precise orders of the monomials.
In particular, the order of $X_1^MX_2^{2^L-2M-1}$ is 16, and hence the class in (\ref{class}) is nonzero
since it has a term $8uX_1^MX_2^{2^L-2M-1}$, and so (\ref{class}) contradicts the hypothesized immersion.

Part b of \ref{h=1} is proved similarly. If $P^{8M+8}$ immerses in $\R^{16M+4}$, then there is an axial map
$$P^{8M+8}\times P^{2^{L+3}-16M-6}@>f>> P^{2^{L+3}-8M-10},$$
and hence, up to odd multiples, 
\begin{eqnarray}\label{eq2}&&2^{2-\nu(M+1)}X_1^{M+1}X_2^{2^L-2M-2}+2^2X_1^MX_2^{2^L-2M-1}\\
&=&0\in\nonumber \tmf^{2^{L+3}-8M-8}(P^{8M+8}\w P^{2^{L+3}-16M-6}),\end{eqnarray}
since $\a(M)=2$. We have $\tmf^{8M+8}(P^{8M+8})\approx \tmf_{-1}(P_{-1})\approx\Z/2$, and
$$\tmf^{2^{L+3}-16M-8}(P^{2^{L+3}-16M-6})\approx \tmf_{-1}(P_{-3})\approx\Z/8.$$
Thus the two monomials in (\ref{eq2}) have order at most 2 and 8, respectively.
On the other hand, the group in (\ref{eq2}) is isomorphic to $\tmf_6(P_{-1}\w P_{-3})$.
A minimal resolution calculation easier than the one in Diagram \ref{diag1} shows that $\tmf_6(P_{-1}\w P_{-3})$
has $\Z/2\oplus\Z/8$ emanating from filtration 0 (and another $\Z/2\oplus\Z/8$ in higher filtration). The monomials
of (\ref{eq2}) are generated in filtration 0, and since the above upper bound for their orders equals the order
of the subgroup generated by filtration-0 classes, we conclude that the orders of the monomials in (\ref{eq2}) are precisely 2 and 8, respectively,
and so the term $4X_1^MX_2^{2^L-2M-1}$ in (\ref{eq2}) is nonzero, contradicting the immersion.

\section{$\tmf$-cohomology of $P^\infty\times P^\infty$}\label{tPPsec}
In this section, we compute $\tmf^*(P^\infty)$ and $\tmf^{8*}(P^\infty\times P^\infty)$ in positive gradings.
These will be used in the next section in studying the axial class in $\tmf$-cohomology.

There is an element $c_4\in\pi_8(\tmf)$ which reduces to $v_1^4\in\pi_8(bo)$; it has Adams filtration 4.
It acts on $\tmf^*(X)$ with degree $-8$. Recall also that $\pi_*(bo)=bo_*$ is as depicted in \ref{bodiag}.
We denote $bo^*=bo_{-*}$.
We use $P_1$ and $P^\infty$ interchangeably.

\begin{thm}\label{tPthm} There is an element $X\in \tmf^8(P_1)$ of Adams filtration $0$, described in \cite{BDM}, such that, in positive dimensions divisible by $8$, $\tmf^*(P_1)$ is isomorphic as an algebra over $\Z_{(2)}[c_4]$ to $\Z_{(2)}[c_4][X]$.
In particular, each $\tmf^{8i}(P_1)$ with $i>0$ is a free abelian group with basis
$\{c_4^jX^{i+j}:\,j\ge0\}$. There is a class $L\in t^0(P_1)$ such that
\begin{itemize} \item $\tmf^0(P_1)$ is a free abelian group with basis $\{L,c_4^jX^j:\,j\ge 1\}$, and
\item $L^2=2L$ and $LX=2X$.
\end{itemize}

\noindent Moreover, in positive dimensions $\tmf^*(P_1)$ is isomorphic as a graded abelian group to $bo^*[X]$,
and is depicted in Diagram \ref{diag3}.\end{thm}
\begin{rmk}{\rm A complete description of $\tmf^*(P_1)$ as a graded abelian group could probably be obtained
using the analysis in the proof which follows, together with the computation of the $E_2$-term of the
ASS converging to $\tmf_*(P_{-1})$, which was given in \cite{ExtA2}. However, this is
quite complicated and unnecessary for this paper, and so will be omitted.}\end{rmk}
\begin{pf}
We begin with the structure as graded abelian group. There are isomorphisms
\begin{equation}\label{duality}\tmf^*(P_1)\approx\lim_{\leftarrow} \tmf^*(P_1^n)\approx\lim_{\leftarrow} \tmf_{-*-1}(P^{-2}_{-n-1})= \tmf_{-*-1}(P^{-2}_{-\infty}).\end{equation}
Since $H^*(\tmf;\zt)\approx A/\!/A_2$, there is a spectral sequence converging to $\tmf_*(X)$ with $E_2(X)=\ext_{A_2}(H^*X,\zt)$. Here $A_2$ is the subalgebra of the mod 2 Steenrod algebra $A$ generated
by $\sq^1$, $\sq^2$, and $\sq^4$. Also $\zt=\bz/2$.

We compute $E_2(P_{-\infty}^{-2})$ from the exact sequence
\begin{equation}\to E_2^{s-1,t}(P_{-1}^\infty)\to E_2^{s,t}(P^{-2}_{-\infty})\to E_2^{s,t}(P_{-\infty}^\infty)@>q_*>> E_2^{s,t}(P_{-1}^{\infty})\to.\label{exact}\end{equation}
It was proved in \cite{LDMA} that 
$$\ext_{A_2}(\bP_{-\infty}^{\infty},\zt)\approx\bigoplus_{i\in\Z}\ext_{A_1}(\Sigma^{8i-1}\zt,\zt).$$
Here we have initiated a notation that $\bP_n^m:= H^*(P_n^m)$.
A complete calculation of $\ext_{A_2}(\bP_{-1}^\infty,\zt)$ was performed in \cite{ExtA2}, but all
we need here are the first few groups.
We can now form a chart for $E_2(P^{-2}_{-\infty})$ from (\ref{exact}), as in Diagram
\ref{diag2}, where $\circ$ indicate elements of $\ext_{A_2}(\bP_{-1}^\infty,\zt)$ suitably positioned, and lines of negative
slope correspond to cases of $q_*\ne0$ in (\ref{exact}).

\begin{center}
\begin{minipage}{6.5in}
\begin{diag}\label{diag2}{$\tmf_*(P_{-\infty}^{-2})$, $-17\le*\le 2$}
\begin{center}
\begin{picture}(440,180)
\def\mp{\multiput}
\def\elt{\circle*{3}}
\def\cir{\circle{3}}
\put(40,15){\line(1,0){370}}
\put(90,0){$-17$}
\put(213,0){$-9$}
\put(333,0){$-1$}
\put(60,60){$\cdots$}
\mp(100,15)(0,15){7}{\elt}
\mp(220,15)(0,15){7}{\elt}
\mp(340,15)(0,15){7}{\elt}
\mp(100,15)(120,0){3}{\vector(0,1){130}}
\mp(100,15)(120,0){3}{\line(1,1){30}}
\mp(115,30)(120,0){3}{\elt}
\mp(130,45)(120,0){3}{\elt}
\mp(103,76)(120,0){3}{\vector(0,1){69}}
\mp(103,76)(120,0){3}{\line(1,1){30}}
\mp(103,76)(120,0){3}{\elt}
\mp(118,91)(120,0){3}{\elt}
\mp(133,106)(120,0){3}{\elt}
\mp(103,91)(120,0){3}{\elt}
\mp(103,106)(120,0){3}{\elt}
\mp(160,60)(120,0){2}{\vector(0,1){85}}
\mp(160,60)(0,15){4}{\elt}
\mp(280,60)(0,15){4}{\elt}
\mp(163,121)(120,0){2}{\vector(0,1){24}}
\mp(163,121)(120,0){2}{\elt}
\put(325,30){\cir}
\mp(355,45)(15,15){3}{\cir}
\mp(385,30)(0,15){3}{\cir}
\mp(357,47)(15,15){2}{\line(1,1){11}}
\mp(385,32)(0,15){3}{\line(0,1){11}}
\put(340,15){\line(-1,1){14}}
\put(385,75){\line(1,-1){15}}
\end{picture}
\end{center}
\end{diag}
\end{minipage}
\end{center}

Dualizing, we obtain Diagram \ref{diag3} for the desired $\tmf^*(P_1^\infty)$.
\begin{center}
\begin{minipage}{6.5in}
\begin{diag}\label{diag3}{$\tmf^*(P_{1}^{\infty})$, $*\ge-2$}
\begin{center}
\begin{picture}(440,180)
\def\mp{\multiput}
\def\elt{\circle*{3}}
\def\cir{\circle{3}}
\put(50,15){\line(1,0){350}}
\put(98,0){$0$}
\put(218,0){$8$}
\put(335,0){$16$}
\mp(100,30)(0,15){6}{\elt}
\mp(220,15)(0,15){7}{\elt}
\mp(340,15)(0,15){7}{\elt}
\mp(85,30)(-15,15){2}{\elt}
\put(85,30){\line(-1,1){15}}
\mp(85,45)(-15,15){2}{\elt}
\put(100,30){\line(-1,1){30}}
\mp(220,15)(120,0){2}{\line(-1,1){30}}
\mp(205,30)(120,0){2}{\elt}
\mp(190,45)(120,0){2}{\elt}
\put(100,30){\vector(0,1){115}}
\mp(220,15)(120,0){2}{\vector(0,1){130}}
\mp(97,76)(120,0){3}{\vector(0,1){69}}
\mp(97,76)(0,15){3}{\elt}
\mp(217,76)(0,15){3}{\elt}
\mp(337,76)(0,15){3}{\elt}
\mp(82,91)(120,0){3}{\elt}
\mp(67,106)(120,0){3}{\elt}
\mp(97,76)(120,0){3}{\line(-1,1){30}}
\mp(160,60)(120,0){2}{\vector(0,1){85}}
\mp(160,60)(0,15){4}{\elt}
\mp(280,60)(0,15){4}{\elt}
\mp(157,121)(120,0){2}{\vector(0,1){24}}
\mp(157,121)(120,0){2}{\elt}
\put(360,60){$\cdots$}
\end{picture}
\end{center}
\end{diag}
\end{minipage}
\end{center}

Naming of the generators $X^i$ is clear since $X$ has filtration 0. The free action
of $c_4$ is also clear. The class $L$ is (up to sign) the composite $P_1@>\lambda>> S^0\to \tmf$, where
$\lambda$ is the well-known Kahn-Priddy map. Thus $L$ is the image of a class $\hat L\in\pi^0(P_1)$.
 Lin's theorem (\cite{Lin}) says that  $\pi^0(P_{1})\approx\zt^\w$, generated by $\hat L$.
Since $\pi^0(P_1)\to ko^0(P_1)$ is an isomorphism, and, since $(1-\xi)^2=2(1-\xi)$ for a generator $(1-\xi)$ of
$ko^0(P_1)$, we obtain $\hat L^2=2\hat L$, and hence also for $L$. We chose the generator to be $(1-\xi)$ rather than $(\xi-1)$ to avoid minus
signs later in the paper.

To prove the claim about $LX$, first note that, by the structure of $\tmf^8(P_1)$, we must have
$LX=p(c_4X)X$ for some polynomial $p$. Multiply both sides by $L$ and apply the result about $L^2$ to get
$2LX=p(c_4X)LX$, hence $2p=p^2$, from which we conclude $p=2$.\end{pf}
 
 In $\tmf^*(P_1\times P_1)$, for $i=1,2$, let $L_i$ and $X_i$ denote the classes $L$ and $X$ in the $i$th factor.
Note that there is an isomorphism as $\tmf_*$-modules, but not as rings, 
$$\tmf^*(P_1\times P_1)\approx \tmf^*(P_1\w P_1)\oplus \tmf^*(P_1\times*)\oplus \tmf^*(*\times P_1).$$
\begin{thm}\label{tPPthm} In positive dimensions divisible by $8$,
$\tmf^*(P_1 \w P_1)$  is isomorphic as a graded abelian group to a free abelian group on monomials $X_1^iX_2^j$ with $i,j>0$ direct sum with a free 
$\Z[c_4]$-module with basis $\{L_1X_2^i,\, X_1^iL_2:\,i\ge 1\}$.
The product and $\Z[c_4]$-module structure is determined from \ref{tPthm} and
$$c_4(X_1X_2)=(c_4X_1)X_2=X_1(c_4X_2)=\sum_{i\ge 0}\g_ic_4^i(L_1X_2^{i+1}+X_1^{i+1}L_2),$$
for certain integers $\g_i$ with $\g_0$ divisible by $8$.
\end{thm}
 
The proof of this theorem involves a number of subsidiary results. They and it occupy the
remainder of this section. 
We will use duality and exact sequences similar
to (\ref{exact}). But to get started, we need $\ext_{A_2}(\bP\ot\bP,\zt)$. Here we have begun to abbreviate
$\bP:=\bP_{-\infty}^\infty$. We begin with a simple lemma. Throughout this section, $x_1$ and $x_2$ denote
nonzero elements coming from the factors in $H^1(RP\times RP;\zt)$. 
\begin{lem}\label{split} $(\cite{Seg})$ There is a split short exact sequence of $A$-modules
$$0\to \zt\ot\bP\to \bP\ot\bP\to(\bP/\zt)\ot\bP\to 0.$$\end{lem}
\begin{pf} The $\zt$ is, of course, the subgroup generated by $x^0$, which is an $A$-submodule.
A splitting morphism $\bP\ot\bP@>g>>\zt\ot\bP$ is defined by $g(x_1^i\ot x_2^j)=x_1^0\ot x_2^{i+j}$.
This is $A$-linear since 
$$g(\sq^k(x_1^i\ot x_2^j))=\sum_\ell\tbinom i\ell\tbinom j{k-\ell}x_1^0\ot x_2^{i+j+k}=\tbinom{i+j}kx_1^0\ot x_2^{i+j+k}=\sq^k
g(x_1^i\ot x_2^j).$$
\end{pf}
The following result is more substantial. We will prove it at the end of this section.
\begin{prop}\label{CB} There is a short exact sequence of $A_2$-modules
$$0\to C\to (\bP/\zt)\ot\bP\to B\to 0,$$
where $C$ has a filtration with
$$F_p(C)/F_{p-1}(C)\approx\Sigma^{8p}A_2/\sq^2,\ p\in\bz,$$
and $B$ has a filtration with 
$$F_p(B)/F_{p-1}(B)\approx\bigoplus_{\bz\text{ copies}}\Sigma^{4p-2}A_2/\sq^1,\ p\in\bz.$$
The generator of $F_p(C)/F_{p-1}(C)$ is $x_1^1x_2^{8p-1}$; a basis over $\zt$ for $C$ is 
$$\{x_1^2x_2^{i+2}+x_1^4x_2^i,x_1^4x_2^i+x_1^8x_2^{i-4},\,i\in\Z\}\cup\{x_1^1x_2^{i-1}+x_1^2x_2^{i-2},\,i\not\equiv0\,(8)\}\cup
\{x_1^1x_2^{8p-1},\, p\in\Z\}.$$
A minimal set of generators as an $A_2$-module for the filtration quotients of $B$ is $\{x_1^{8i-1}x_2^{4j-1}:\,i,j\in\Z\}$.
\end{prop}
\begin{cor} A chart for $\ext_{A_2}^{s,t}(\bP\ot\bP,\zt)$ in $8p-3\le t-s\le 8p+4$ is as suggested in Diagram \ref{diag4}, for all integers $p$. The big batch of towers in each grading $\equiv2\,(4)$ represents an infinite family of towers. The pattern of the other classes is repeated with vertical period 4. Thus, for example, in $8p-1$ there is
an infinite tower emanating from filtration $4i$ for each $i\ge0$.
\label{batch}\end{cor}
\begin{center}
\begin{minipage}{6.5in}
\begin{diag}\label{diag4}{$\ext_{A_2}^{s,t}(\bP\ot\bP,\zt)$ in $8p-3\le t-s\le 8p+4$}
\begin{center}
\begin{picture}(440,280)
\def\mp{\multiput}
\def\elt{\circle*{3}}
\put(55,30){\line(1,0){245}}
\put(55,15){$8p+$}
\put(94,15){$-2$}
\put(157,15){$0$}
\put(217,15){$2$}
\put(277,15){$4$}
\mp(95,30)(1,0){11}{\vector(0,1){240}}
\mp(215,30)(1,0){11}{\vector(0,1){240}}
\put(128,30){\vector(0,1){240}}
\put(152,120){\vector(0,1){150}}
\put(131,150){\vector(0,1){120}}
\put(128,30){\line(1,1){60}}
\put(152,120){\line(1,1){60}}
\put(131,150){\line(1,1){60}}
\mp(128,30)(30,30){3}{\elt}
\mp(152,120)(30,30){3}{\elt}
\mp(131,150)(30,30){3}{\elt}
\mp(182,30)(30,30){2}{\elt}
\put(182,30){\line(1,1){30}}
\put(280,90){\vector(0,1){180}}
\put(250,120){\vector(0,1){150}}
\put(283,210){\vector(0,1){60}}
\put(253,240){\vector(0,1){30}}
\put(280,90){\elt}
\put(250,120){\elt}
\put(283,210){\elt}
\put(253,240){\elt}
\end{picture}
\end{center}
\end{diag}
\end{minipage}
\end{center}
 
\begin{pf*}{Proof of Corollary \ref{batch}} We first note that $\ext_{A_2}(\bP,\zt)$ is identical to the left portion
of Diagram \ref{diag2} extended periodically in both directions. Also, $\ext_{A_2}(A_2/\sq^1,\zt)\approx\ext_{A_0}(\zt,\zt)$
is just an infinite tower, and $$\ext_{A_2}(A_2/\sq^2,\zt)\approx\ext_{A_1}(A_1/\sq^2,\zt)$$
is given as in Diagram \ref{diag5}. We will show at the end of this proof that
\begin{equation}\label{Ceq}\ext_{A_2}(C,\zt)\approx\bigoplus_{p\in\Z}\ext_{A_2}(\Sigma^{8p}A_2/\sq^2,\zt)\end{equation}
and similarly  $$\ext_{A_2}(B,\zt)\approx\bigoplus_p\bigoplus_{\Z}\ext_{A_2}(\Sigma^{4p-2}A_2/\sq^1,\zt).$$ These would follow by induction on $p$ once you get started,
but since $p$ ranges over all integers, that is not automatic.

Thus $\ext_{A_2}(\bP\ot\bP,\zt)$ is formed from
$$\ext_{A_2}(\bP,\zt)\oplus\bigoplus\ext_{A_2}(\Sigma^{8p}A_2/\sq^2,\zt)\oplus\bigoplus\ext_{A_2}(\Sigma^{4p-2}A_2/\sq^1,\zt),$$
using the sequences in \ref{split} and \ref{CB}. The Ext sequence of \ref{split} must split,
and there are no possible boundary morphisms in the Ext sequence of \ref{CB}, yielding the claim of the corollary.

To prove (\ref{Ceq}), let $(s,t)$ be given, and choose $p_0$ so that $8p_0<t-23s+2$. Since the highest degree element in $A_2$ is in degree 23,
$\ext_{A_2}^{s,t}(F_{p_0}(C),\zt)=0$. Actually a much sharper lower vanishing line
can be established, but this is good enough for our purposes. Thus, for this $(s,t)$,
\begin{equation}\label{p11}\ext_{A_2}^{s,t}(F_{p_1}(C),\zt)\approx\bigoplus_{p\le p_1}\ext_{A_2}^{s,t}(\Sigma^{8p-2}A_2/\sq^2)\end{equation}
for $p_1\le p_0$, as both are 0.
Let $p_1$ be minimal such that (\ref{p11}) does not hold. Then comparison of exact sequences implies that
$$\ext_{A_2}^{s-1,t}(F_{p_1-1}(C),\zt)\to \ext_{A_2}^{s,t}(F_{p_1}(C)/F_{p_1-1}(C),\zt)$$
must be nonzero. But one or the other of these groups is always 0,\footnote{Actually this is not quite true; for one
family of elements we need to use $h_0$-naturality.}
 as both charts $\ext_{A_2}^{*,*}(F_{p_1-1}(C),\zt)$
and $\ext_{A_2}^{*,*}(F_{p_1}(C)/F_{p_1-1}(C),\zt)$ are copies of Diagram \ref{diag5} displaced by 4 vertical
units from one another. Thus (\ref{p11}) is true for all $p_1$, and hence (\ref{Ceq}) holds.
A similar proof works when $C$ is replaced by $B$.
\end{pf*} 
\begin{center}
\begin{minipage}{6.5in}
\begin{diag}\label{diag5}{$\ext_{A_2}(A_2/\sq^2,\zt)$}
\begin{center}
\begin{picture}(440,100)
\def\mp{\multiput}
\def\elt{\circle*{3}}
\put(158,10){$0$}
\put(150,20){\line(1,0){100}}
\mp(160,20)(10,10){2}{\elt}
\put(160,20){\line(1,1){10}}
\put(190,40){\vector(0,1){45}}
\mp(190,40)(0,10){3}{\elt}
\put(230,50){\vector(0,1){35}}
\mp(230,50)(0,10){3}{\elt}
\put(230,50){\line(1,1){20}}
\mp(240,60)(10,10){2}{\elt}
\put(265,85){$\cdots$}
\end{picture}
\end{center}
\end{diag}
\end{minipage}
\end{center}

Now we can prove a result which will, after dualizing, yield Theorem \ref{tPPthm}. The groups $\ext_{A_1}(\zt,\zt)$
to which it alludes are depicted in \ref{bodiag}. The content of this result is pictured in Diagram
\ref{diag8}.
\begin{prop}\label{P-2} In dimensions $t-s\equiv 2\mod4$ with $t-s\le-10$, $\ext_{A_2}(\bP^{-2}_{-\infty}\ot
\bP^{-2}_{-\infty},\zt)$ consists of $i$ infinite towers emanating from filtration $0$ in dimensions $-8i-6$ and $-8i-10$,
together with the relevant portion of two copies of $\ext_{A_1}(\zt,\zt)$ beginning in filtration 1 in each dimension $-8i-2$.
The generators of the towers in $-8i-10$ correspond to cohomology classes $x_1^{-9}x_2^{-8i-1},\ldots,x_1^{-8i-1}x_2^{-9}$.
The generators of the two copies of $\ext_{A_1}(\zt,\zt)$ in $-8i-2$ arise from $h_0$ times classes corresponding
to $x_1^{-1}x_2^{8i-1}$ and $x_1^{-8i-1}x_2^{-1}$.\end{prop}
\begin{pf} Using exact sequences like (\ref{exact}) on each factor, we build $\ext_{A_2}^{*,*}(\bP^{-2}_{-\infty}\ot
\bP^{-2}_{-\infty},\zt)$ from $\bA:=\ext_{A_2}^{*,*}(\bP\ot\bP,\zt)$, $\bB:=\ext_{A_2}^{*-1,*}(\bP^{\infty}_{-1}\ot\bP,\zt)$, $\bC:=\ext_{A_2}^{*-1,*}(\bP\ot\bP_{-1}^{\infty},\zt)$,
and $\bD:=\ext_{A_2}^{*-2,*}(P_{-1}^{\infty}\ot P_{-1}^{\infty},\zt)$, with possible $d_1$-differential
from $\bA$ and into $\bD$. In the range of concern, $t-s\le-9$, the $\bD$-part will not be present, and the
part of Diagram \ref{diag4} in dimension $\not\equiv2$ mod 4 will not be involved in $d_1$. Using \cite{LDMA} for $\bB$ and $\bC$,
the relevant part, namely the portion of $\bA$ in dimension $\equiv2$ mod 4, together with $\bB$ and $\bC$, is pictured
in Diagram \ref{diag6}.

\begin{center}
\begin{minipage}{6.5in}
\begin{diag}\label{diag6}{Portion of $\bA+\bB+\bC$}
\begin{center}
\begin{picture}(440,220)
\def\mp{\multiput}
\def\elt{\circle*{3}}
\put(90,20){\line(1,0){190}}
\mp(120,20)(1,0){10}{\vector(0,1){200}}
\mp(200,20)(1,0){10}{\vector(0,1){200}}
\mp(280,20)(1,0){10}{\vector(0,1){200}}
\put(120,5){$-2$}
\put(203,5){$2$}
\put(283,5){$6$}
\put(60,5){$8p+$}
\mp(100,40)(2,0){2}{\elt}
\mp(260,40)(2,0){2}{\elt}
\mp(104,40)(0,20){5}{\elt}
\mp(264,40)(0,20){5}{\elt}
\mp(106,40)(0,20){5}{\elt}
\mp(266,40)(0,20){5}{\elt}
\mp(104,40)(2,0){2}{\line(0,1){80}}
\mp(264,40)(2,0){2}{\line(0,1){80}}
\mp(108,40)(0,20){9}{\elt}
\mp(268,40)(0,20){9}{\elt}
\mp(110,40)(0,20){9}{\elt}
\mp(270,40)(0,20){9}{\elt}
\mp(108,40)(2,0){2}{\line(0,1){160}}
\mp(268,40)(2,0){2}{\line(0,1){160}}
\mp(142,60)(2,0){2}{\elt}
\mp(161,79)(2,0){2}{\elt}
\mp(180,40)(0,20){4}{\elt}
\mp(182,40)(0,20){4}{\elt}
\mp(180,40)(2,0){2}{\line(0,1){60}}
\mp(142,60)(2,0){2}{\line(1,1){38}}
\mp(142,140)(2,0){2}{\elt}
\mp(163,161)(2,0){2}{\elt}
\mp(184,40)(0,20){8}{\elt}
\mp(186,40)(0,20){8}{\elt}
\mp(184,40)(2,0){2}{\line(0,1){140}}
\mp(142,140)(2,0){2}{\line(1,1){42}}
\end{picture}
\end{center}
\end{diag}
\end{minipage}
\end{center}

In dimension $8p-2$, the towers in $\bA$ arise from all cohomology classes $x_1^{-8i-1}x_2^{-8j-1}$ with $i+j=-p$,
while in dimension $8p+2$, they arise from $x_1^{8i-1}x_2^{8j+3}\sim x_1^{8i+3}x_2^{8j-1}$. The finite towers in $\bB$ arise from
$x_1^{4i-1}x_2^{8j-1}$ with $i\ge0$, and those from $\bC$ from $x_1^{8i-1}x_2^{4j-1}$ with $j\ge0$. The homomorphism
$$\ext_{A_2}^0(\bP\ot \bP,\zt)\to\ext_{A_2}^0(\bP_{-1}^{\infty}\ot\bP,\zt)\oplus\ext_{A_2}^0(\bP\ot\bP_{-1}^{\infty},\zt),$$
which is equivalent to the $d_1$-differential mentioned above, sends classes to those with the same name.
In dimension $\le-10$, this is surjective, with kernel spanned by classes with both components $<-1$. In dimension
$-8i-6$ and $-8i-10$, there will be $i$ such classes. We illustrate by listing the classes in the first few gradings:
\begin{eqnarray*}-14&:&x_1^{-9}x_2^{-5}\sim x_1^{-5}x_2^{-9}\\
-18&:&x_1^{-9}x_2^{-9}\\
-22&:&x_1^{-17}x_2^{-5}\sim x_1^{-13}x_2^{-9},\ x_1^{-9}x_2^{-13}\sim x_1^{-5}x_2^{-17}\\
-26&:&x_1^{-17}x_2^{-9},\ x_1^{-9}x_2^{-17}.\end{eqnarray*}
These kernel classes yield infinite towers emanating from filtration 0. 

For each $p<0$, the towers arising from $x_1^{4j-1}x_2^{8p-1}$, $j\ge0$, in $\bA$ combine with those in the $p$-summand of
$$\bB\approx\bigoplus_{p\in\Z}\ext_{A_1}(\Sigma^{8p-1}P_{-1}^{\infty},\zt)$$
as in Diagram \ref{diag7} to yield one of the copies of $\ext_{A_1}(\zt,\zt)$ arising from filtration 1.
An identical picture results when the factors are reversed.
\end{pf}

\begin{center}
\begin{minipage}{6.5in}
\begin{diag}\label{diag7}{Part of $\ext_{A_2}(\bP^{-2}_{-\infty}\ot\bP^{-2}_{-\infty},\zt)$}
\begin{center}
\begin{picture}(440,140)
\def\mp{\multiput}
\def\elt{\circle*{3}}
\put(0,0){\line(1,0){165}}
\put(270,0){\line(1,0){170}}
\mp(15,0)(0,15){8}{\elt}
\mp(15,0)(60,0){3}{\vector(0,1){120}}
\mp(75,0)(0,15){8}{\elt}
\mp(135,0)(0,15){8}{\elt}
\put(0,15){\elt}
\put(0,15){\line(1,-1){15}}
\mp(75,0)(0,15){4}{\line(-1,1){15}}
\mp(60,15)(0,15){4}{\elt}
\put(60,15){\line(0,1){45}}
\mp(30,30)(15,15){2}{\elt}
\put(30,30){\line(1,1){30}}
\mp(135,0)(0,15){5}{\line(-1,1){15}}
\mp(120,15)(0,15){5}{\elt}
\put(120,15){\line(0,1){60}}
\put(150,90){\line(1,1){15}}
\mp(150,90)(15,15){2}{\elt}
\put(220,60){$\Longrightarrow$}
\mp(280,15)(0,15){7}{\elt}
\put(280,15){\vector(0,1){105}}
\mp(295,30)(15,15){2}{\elt}
\put(280,15){\line(1,1){30}}
\put(340,60){\vector(0,1){60}}
\mp(340,60)(0,15){4}{\elt}
\put(400,75){\vector(0,1){45}}
\mp(400,75)(0,15){3}{\elt}
\put(400,75){\line(1,1){30}}
\mp(415,90)(15,15){2}{\elt}
\end{picture}
\end{center}
\end{diag}
\end{minipage}
\end{center}

\medskip
Putting things together, we obtain that in dimensions less than $-8$,
$\ext_{A_2}(\bP^{-2}_{-\infty}\ot\bP^{-2}_{-\infty},\zt)$ consists of
a chart described in Proposition \ref{P-2} and partially illustrated in
Diagram \ref{diag8} together with the classes in Diagram \ref{diag4}
which are not part of the infinite sums of towers in dimension $\equiv2$ mod 4.

\begin{center}
\begin{minipage}{6.5in}
\begin{diag}\label{diag8}{Illustration of Proposition \ref{P-2}}
\begin{center}
\begin{picture}(450,250)
\def\mp{\multiput}
\def\elt{\circle*{3}}
\put(0,25){\line(1,0){450}}
\put(-7,10){$-26$}
\put(193,10){$-18$}
\put(392,10){$-10$}
\mp(0,25)(2,0){2}{\vector(0,1){225}}
\mp(4,50)(2,0){2}{\vector(0,1){200}}
\mp(4,50)(2,0){2}{\line(1,1){50}}
\mp(8,150)(2,0){2}{\vector(0,1){100}}
\mp(8,150)(2,0){2}{\line(1,1){50}}
\mp(100,25)(2,0){2}{\vector(0,1){225}}
\mp(104,125)(2,0){2}{\vector(0,1){125}}
\mp(108,225)(2,0){2}{\vector(0,1){25}}
\mp(200,25)(100,0){2}{\vector(0,1){225}}
\mp(202,50)(2,0){2}{\vector(0,1){200}}
\mp(202,50)(2,0){2}{\line(1,1){50}}
\mp(206,150)(2,0){2}{\vector(0,1){100}}
\mp(206,150)(2,0){2}{\line(1,1){50}}
\mp(302,125)(2,0){2}{\vector(0,1){125}}
\mp(306,225)(2,0){2}{\vector(0,1){25}}
\mp(400,50)(2,0){2}{\vector(0,1){200}}
\mp(400,50)(2,0){2}{\line(1,1){50}}
\mp(404,150)(2,0){2}{\vector(0,1){100}}
\mp(404,150)(2,0){2}{\line(1,1){50}}
\end{picture}
\end{center}
\end{diag}
\end{minipage}
\end{center}

The only possible differentials in the Adams spectral sequence of
$P^{-2}_{-\infty}\w P^{-2}_{-\infty}\w\tmf$ involving the classes in
dimensions $8p-2$ with $p<0$ are from the towers in $8p-1$ in
Diagram \ref{diag4}, but these differentials are shown to be 0 as in \cite[p.54]{BDM}.
Similarly to (\ref{duality}), we have
$$\tmf^*(P_1\w P_1)\approx \tmf_{-*-2}(P^{-2}_{-\infty}\w P^{-2}_{-\infty}),$$
and so we obtain a turned-around version of Diagram \ref{diag8}, of the same general sort
as Diagram \ref{diag3}, as a depiction
of a relevant portion of $\tmf^*(P_1\w P_1)$, with the labeled columns in Diagram \ref{diag8}
corresponding to cohomology gradings 24, 16, and 8.

The classes $X_1^iX_2^j$ described in Theorem \ref{tPPthm} are detected
by the S-duals of the classes from which the filtration-0 towers in
dimensions $8p-2$ in Diagram \ref{diag8} arise, and so they can be chosen to be the
corresponding elements of $\tmf^{8*}(P_1\w P_1)$. Similarly the classes
$L_1X_2^i$ and $X_1^iL_2$ have Adams filtration 1, and so one would
anticipate that they represent the duals of the generators of the two towers in dimension
$8p-2$ with $p<0$ in Diagram \ref{diag8}. This seems a bit harder to
prove using the Adams spectral sequence; however, the Atiyah-Hirzebruch
spectral sequence shows this quite clearly. The class $X_1^i$ is detected by
$H^{8i}(P_1;\pi_0(\tmf))$, while $L$ is detected by $H^1(P_1;\pi_1(\tmf))$.
Under the pairing, their product is detected in $H^{8i+1}(P_1;\pi_1(\tmf))$,
clearly of Adams filtration 1.

The last part of Theorem \ref{tPPthm} deals with the action of $c_4$ on the
monomials $X_1^iX_2^j$. Since $\tmf$ is a commutative ring spectrum, $\tmf^*(P_1\w P_1)$
is a graded commutative algebra over $\tmf_*$. The action $c_4(X_1X_2)$ must be of
the form $\sum_{i\ge0}\g_ic_4^i(L_1X_2^i+X_1^iL_2)$ as these are the only elements
in $\tmf^8(P_1\w P_1)$, and the class must be invariant under reversing factors.
The divisibility of $\g_0$ by 8 follows since $c_4$ has Adams filtration 4.

Having just completed the proof of Theorem \ref{tPPthm}, we conclude this section with
the postponed proof of Proposition \ref{CB}.

\begin{pf*}{Proof of Proposition \ref{CB}} Let $C$ denote the $A_2$-submodule of $(\bP/\zt)\ot\bP$
generated by all $x_1^{1}x_2^{8p-1}$, $p\in\Z$.
Note that $\sq^2(x_1^{1}x_2^{8p-1})=\sq^4\sq^6(x_1^1x_2^{8p-9})$.
Thus a basis of $A_2/\sq^2$ acting on all $x_1^1x_2^{8p-1}$ spans $C$.
The 24 elements in a basis of $A/\sq^2$ acting on $x_1^1x_2^7$ yield $x_1^1x_2^7$, $x_1^1x_2^8+x_1^2x_2^7$, $x_1^2x_2^9+x_1^4x_2^7$,
$x_1^1x_2^{11}+x_1^2x_2^{10}$, $x_1^1x_2^{12}+x_1^2x_2^{11}$, $x_1^1x_2^{13}+x_1^2x_2^{12}$, $x_1^1x_2^{14}+x_1^2x_2^{13}$, $x_1^2x_2^{13}+x_1^4x_2^{11}$, $x_1^4x_2^{11}+x_1^8x_2^7$,
$x_1^2x_2^{14}+x_1^4x_2^{12}$, $x_1^1x_2^{16}+x_1^4x_2^{13}$, $x_1^1x_2^{17}+x_1^2x_2^{16}$, $x_1^2x_2^{16}+x_1^4x_2^{14}$, $x_1^1x_2^{18}+x_1^2x_2^{17}$, $x_1^2x_2^{17}+x_1^8x_2^{11}$,
$x_1^2x_2^{18}+x_1^8x_2^{12}$, $x_1^1x_2^{20}+x_1^4x_2^{17}$, $x_1^4x_2^{17}+x_1^8x_2^{13}$, $x_1^2x_2^{20}+x_1^4x_2^{18}$, $x_1^4x_2^{18}+x_1^8x_2^{14}$, $x_1^4x_2^{20}+x_1^8x_2^{16}$,
$x_1^1x_2^{24}+x_1^8x_2^{17}$, $x_1^2x_2^{24}+x_1^8x_2^{18}$, and $x_1^4x_2^{24}+x_1^8x_2^{20}$. These classes with second components shifted
by all multiples of 8 exactly comprise the basis for $C$ described in the proposition.

The procedure to establish the structure of $B=((\bP/\zt)\ot\bP)/C$ is similar but more
elaborate. For the 32 elements $\theta$ in a basis of $A_2/\sq^1$, we list $\theta(x_1^{-1}x_2^{-1})$
and $\theta(x_1^{-1}x_2^3)$. Then we show that these, with each component allowed to vary by multiples
of 8, together with $C$, fill out all of $(\bP/\zt)\ot\bP$. 

It is convenient to let $\bQ$ denote the quotient of $(\bP/\zt)\ot\bP$ by $C$ and all elements $\theta(x_1^{8i-1}x_2^{8j-1})$
and $\theta(x_1^{8i-1}x_2^{8j+3})$. We will show $\bQ=0$. This will complete the proof of Proposition \ref{CB}, implying in particular that
$\sq^1(x_1^{8i-1}x_2^{8j-1})$ and $\sq^1(x_1^{8i-1}x_2^{8j+3})$ are decomposable over $A_2$.

A separate calculation is performed
for each mod 8 value of the degree. Here we use repeatedly that the $A_2$-action on $x^i$ depends
only on $i$ mod 8.
We illustrate with the case in which degree $\equiv0$ mod 8. The other 7 congruences are handled
similarly, although some are a bit more complicated.

A basis of $A_2/\sq^1$ in degree $\equiv2$ mod 8 acting on $x_1^{-1}x_2^{-1}$ yields the following elements:
$x_1^{-1}x_2^1+x_1^0x_2^0+x_1^1x_2^{-1}$, $x_1^2x_2^6+x_1^6x_2^2$, $x_1^{-1}x_2^9+x_1^3x_2^5+x_1^4x_2^4+x_1^5x_2^3+x_1^9x_2^{-1}$, and $x_1^4x_2^{12}+x_1^{12}x_2^4$.
A basis of $A_2/\sq^1$ in degree $\equiv6$ mod 8 acting on $x_1^{-1}x_2^3$ yields the following elements:
$x_1^2x_2^6+x_1^3x_2^5+x_1^4x_2^4+x_1^5x_2^3$, $x_1^{-1}x_2^9+x_1^2x_2^6+x_1^5x_2^3$, $x_1^4x_2^{12}+x_1^6x_2^{10}+x_1^{10}x_2^6+x_1^{12}x_2^4$, and $x_1^8x_2^{16}+x_1^{16}x_2^8$.
Because we allow both components to vary by multiples of 8, we will list just the first component of the ordered
pairs. These are considered as relations in $\bold{Q}$.
Thus the relation $R_1$ below really means that all $x_1^{8i-1}x_2^{8j+1}+x_1^{8i}x_2^{8j}+x_1^{8i+1}x_2^{8j-1}$ become 0 in $\bold Q$.
\begin{eqnarray*}R_1&:&X_{-1}+X_0+X_1,\\
R_2&:&X_2+X_6,\\
R_3&:&X_{-1}+X_3+X_4+X_5+X_9,\\
R_4&:& X_4+X_{12},\\
R_5&:&X_2+X_3+X_4+X_5,\\
R_6&:&X_{-1}+X_2+X_5,\\
R_7&:&X_4+X_6+X_{10}+X_{12},\\
R_8&:&X_8+X_{16}.\end{eqnarray*}

We will use these relations to show that all classes (in degree $\equiv0$ mod 8) are 0 in $\bold Q$.
First, $R_8$ implies that all classes $X_{8i}$ are congruent to one another. Since $X_0$ is 0 in the
quotient due to $\bP/\zt$, we conclude that all classes $X_{8i}$ are 0 in $\bold Q$. Next, $R_4$ implies
that all $X_{8i+4}$ are congruent to one another. Since $X_4+X_8\in C$, and we have just shown that $X_8\equiv0$
in $\bold Q$, we deduce that all $X_{8i+4}$ are 0 in $\bold Q$. Now we use $R_2+R_7$ to see that all $X_{8i+2}+X_{8i+4}$
are congruent to one another, then that $X_2+X_4\in C$ to deduce all $X_{8i+2}+X_{8i+4}\equiv 0$, and finally
the result of the previous sentence to conclude all $X_{8i+2}\equiv0$. Then $R_2$ implies all $X_{8i+6}\equiv0$.
Now $R_1+R_3+R_5$, together with relations previously obtained, implies all $X_{8i+1}$ are congruent to one another,
and since $X_1\in C$, we conclude all $X_{8i+1}\equiv0$. Finally $R_1$ implies $X_{8i-1}\equiv0$, $R_6$ implies
$X_{8i+5}\equiv0$, and then $R_3$ implies $X_{8i+3}\equiv0$.
\end{pf*}

\section{Careful treatment of axial class}\label{axialsec}
In this section, we fill the gap in the proof in \cite{BDM} of its Theorem 1.1
by careful consideration of the possible ``other terms" in the axial class
discussed in the Introduction. We show that, at least as far as the monomials $cX_1^iX_2^j$ in its powers are concerned, the axial class equals
$\bu(X_1+X_2)$, where $\bu$ is a unit in $\tmf^0(RP^\infty\times RP^\infty)$.
Thus the $\ell$th power of the axial class is nonzero in $\tmf^{8\ell}(RP^n\times RP^m)$
if and only if $(X_1+X_2)^\ell$ is nonzero there, and the latter is the condition
which yielded the nonimmersions of \cite[1.1]{BDM}. Thus we have a complete proof of \cite[1.1]{BDM}.

If $P^n\times P^m@>f>>P^{m+k}$ is an axial map, then
there  is a commutative diagram
$$\begin{CD} P^n\times P^m@>f>> P^{m+k}\\
@VVV @VVV\\
P^\infty\times P^\infty@>g>>P^\infty,\end{CD}$$
where $g$ is the standard multiplication of $P^\infty$, since $P^\infty=K(\zt,1)$.
Since $X\in \tmf^8(P^{m+k})$ has been chosen to  extend over $P^\infty$, we obtain that
$f^*(X)$ is the restriction of $g^*(X)$. By Theorem \ref{tPPthm} and the symmetry of $g$,
we must have
\begin{equation}\label{axialsum}g^*(X)=X_1+X_2+\sum_{i\ge0}\k_ic_4^i(L_1X_2^{i+1}+X_1^{i+1}L_2),\end{equation}
for some integers $\k_i$. This is what we call the ``axial class."
Then $g^*(X^\ell)$ equals the $\ell$th power of (\ref{axialsum}). Using the formulas
for $L_i^2$, $L_iX_i$, and $c_4(X_1X_2)$ in \ref{tPthm} and \ref{tPPthm} and the binomial theorem,
this $\ell$th power can be written in terms of the basis described in \ref{tPPthm}.
If some $\k_i$'s are nonzero, the coefficients of $X_1^iX_2^{\ell-i}$ in $g^*(X^\ell)$ will not equal
$\binom{\ell}i$, as was claimed in \cite{BDM}.
We will study this possible deviation carefully.

One simplification is to treat $L_1$ and $L_2$ as being just 2.
Note that $L_i$ acts like 2 when multiplying by $X_i$, and if, for example, $L_1$ is present
without $X_1$, then the terms $c_4^i L_1 X_2^j$ cannot cancel our
$X_1^kX_2^{\ell}$-classes because both are separate parts of the basis.  You have
to carry the terms along, because they might get multiplied by an $X_1$,
and then it is as if $L_1=2$.  We will incorporate this important simplification throughout
the remainder of this section.

For example, one easily checks that, using $L_1^2=2L_1$ and $L_1X_1=2X_1$, we obtain
$$(X_1+X_2+L_1X_2)^4=(X_1+3X_2)^4-80X_2^4+40L_1X_2^4.$$
The exponent of 2 in each monomial of $(X_1+3X_2)^4-80X_2^4$ is the same as that in $(X_1+X_2)^4$,
and $L_1X_2^4$ is a separate basis element.

With this simplification, the axial class in (\ref{axialsum}) becomes \begin{equation}\label{newaxial}X_1+X_2+2\sum_{i>0}\kappa_i c_4^i(X_1^{i+1}+X_2^{i+1})\end{equation} for some
integers $\kappa_i$.  There was another term $2\kappa_0(X_1+X_2)$, but 
it can be incorporated into the leading $(X_1+X_2)$.  The odd multiple that it can create is not
important.  

{}From Theorem \ref{tPPthm}, we have 
\begin{equation}\label{c41}c_4(X_1X_2)= 16(X_1+X_2)+2\sum_{k>0}\gamma_k c_4^k(X_1^{k+1}+X_2^{k+1}),\end{equation} for some
integers $\g_k$. The 16 comes from $\g_0=8$ and $L_i=2$. Actually we don't really know that
$\g_0=8$, even just up to multiplication by a unit, but it is divisible by 8 and the possibility 
of equality must
be allowed for. This gives
\begin{equation}\label{cxx}c_4(X_1^{i+1}X_2^{j+1})=16(X_1^{i+1}X_2^j+X_1^iX_2^{j+1})+2\sum_{k>0}\g_k c_4^k(X_1^{i+k+1}X_2^j+X_1^iX_2^{j+k+1}).\end{equation}
Here we use that in a graded $\tmf_*$-algebra $\tmf^*(X)$ with even-degree elements,
$c(xy)=cx\cdot y$, for $c\in \tmf_*$ and $x,y\in \tmf^*(X)$.

There is an iterative nature to the action of $c_4$ in (\ref{cxx}), but the leading coefficient 16
enables us to keep track of 2-exponents of leading terms in the iteration. (As observed above, the leading
coefficient might be an even multiple of 16, which would make the terms even more highly 2-divisible. We assume the worst, that it equals 16.) We obtain the following
key result about the action of $c_4$ on monomials in $X_1$ and $X_2$.
\begin{thm} \label{c4thm} There are $2$-adic integers $A_i$ such that
$$c_4=\sum_{i\ge0}2^{4+i}A_i\biggl(\frac 1{X_1}\bigl(\frac{X_2}{X_1}\bigr)^i+\frac 1{X_2}\bigl(\frac{X_1}{X_2}\bigr)^i\biggr).$$
\end{thm}
\begin{rmk}{\rm This formula will be evaluated on (i.e. multiplied by) monomials $X_1^kX_2^\ell$. One might worry that
the negative powers of $X_1$ or $X_2$ in \ref{c4thm} will cause nonsensical negative powers in $c_4X_1^kX_2^\ell$.
This will, in fact, not occur because the monomials on which we act always have total degree greater than the
dimension of either factor. Thus if, after multiplication by $c_4$, a term with negative exponent of $X_i$ appears,
then the accompanying $X_{3-i}^j$-term will be 0 for dimensional reasons.}\end{rmk}
\begin{pf*}{Proof of Theorem \ref{c4thm}}
The defining equation (\ref{c41}) may be written as, with $\theta=c_4\sqrt{X_1X_2}$ and $z=\sqrt{X_1/X_2}$,
\begin{equation}\label{sqrt}\theta=16(z+z^{-1})+\sum_{i>0}2\g_i\theta^i(z^{i+1}+z^{-(i+1)}).\end{equation}
Let $p_i=z^i+z^{-i}$. We will show that \begin{equation}\label{theteq}\theta=\sum_{i\ge0}2^{4+i}A_ip_{2i+1}\end{equation} 
for certain 2-adic integers $A_i$, which interprets back to the claim
of \ref{c4thm}.

 Note that $p_ip_j=p_{i+j}+p_{|i-j|}$, and hence
$$p_1^{e_1}\cdots p_k^{e_k}=p_{\Sigma ie_i}+\cal{L},$$
where $\cal L$ is a sum of integer multiples of $p_j$ with $j<\sum ie_i$ and $j\equiv\sum ie_i$ mod 2.
We will ignore for awhile the coefficients $\g_i$ which occur in (\ref{sqrt}). This is allowable if we agree that 
when collecting terms, we only make crude estimates about their 2-divisibility. We have
\begin{eqnarray*}\theta&=&16p_1+2\theta p_2+2\theta^2p_3+2\theta^3p_4+\cdots\\
&=&16p_1+2p_2(16p_1+2p_2(16p_1+\cdots)+2p_3(16p_1+\cdots)^2+\cdots)\\
&&+2p_3(16p_1+2p_2(16p_1+\cdots)+\cdots)^2+\cdots.\end{eqnarray*}
Note that the only terms that actually get evaluated must end with a $16p_1$ factor.

Now let $T_1=16p_1$ and, for $i\ge2$, let $T_i=2\theta^{i-1}p_i$. Each term in the expansion of $\theta$
involves a sequence of choices. First choose $T_i$ for some $i\ge1$, and then if $i>1$ choose $(i-1)$ factors
$T_j$, one from each factor of $\theta^{i-1}$. For each of these $T_j$ with $j>1$, choose $j-1$ additional factors,
and continue this procedure.
This builds a tree, and we don't get an explicit product term until every branch ends with $T_1$. Each selected 
factor $T_j$ with $j>1$ contributes a factor $2p_j$. There will also be binomial coefficients and the omitted
$\g_i$'s occurring as additional factors.

For example, Diagram \ref{Tdiag} illustrates the choices leading to one term in the expansion of $\theta$.
This yields the term $2p_2\cdot2p_4\cdot16p_1\cdot2p_2\cdot16p_1\cdot2p_3\cdot16p_1\cdot2p_2\cdot16p_1$,
which equals $2^{21}(p_{17}+\cal{L})$, where $\cal L$ is a sum of $p_i$ with $i<17$ and $i$ odd. By induction, one sees
in general that the sum of the subscripts emanating from any node, including the subscript of the node itself, is odd.

\begin{center}
\begin{minipage}{6.5in}
\begin{diag}\label{Tdiag}{A possible choice of terms}
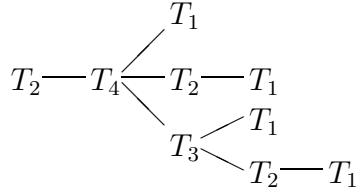
\begin{center}
\begin{picture}(440,80)
\def\mp{\multiput}
\put(150,45){$T_2$}
\mp(162,50)(30,0){3}{\line(1,0){16}}
\put(180,45){$T_4$}
\put(210,45){$T_2$}
\put(240,45){$T_1$}
\put(192,52){\line(1,1){16}}
\put(210,70){$T_1$}
\put(192,48){\line(1,-1){16}}
\put(210,20){$T_3$}
\put(222,27){\line(2,1){16}}
\put(222,23){\line(2,-1){16}}
\put(240,30){$T_1$}
\put(240,10){$T_2$}
\put(252,15){\line(1,0){16}}
\put(270,10){$T_1$}
\end{picture}
\end{center}
\end{diag}
\end{minipage}
\end{center}

The important terms are those in which $T_2$ is chosen $k$ times ($k\ge0$) and then $T_1$ is chosen. These give 
$(2p_2)^kp_1$ with no binomial coefficient. This term is $2^{k+4}(p_{2k+1}+\cal L)$. Note that a term $2^{k+4}p_{2i+1}$
with $i<k$ obtained from $\cal L$ will be more 2-divisible than the $2^{i+4}p_{2i+1}$ term that was previously obtained.
Thus it may be incorporated into the coefficient of that term.

All other terms will be more highly 2-divisible than these. For example, the first would arise from choosing
$T_3$ then two copies of $T_1$. This would give $2p_3\cdot 2^4p_1\cdot 2^4p_1=2^9p_5+\cal L$, and the $2^9p_5$ can
be combined with the $2^6p_5$ obtained from choosing $T_2$ then $T_2$ then $T_1$. Incorporating $\g_i$'s may make
terms even more divisible, but the claim of (\ref{theteq}) is only that $p_{2i+1}$ occurs with coefficient divisible by $2^{4+i}$.
\end{pf*}

Now we incorporate \ref{c4thm} into (\ref{newaxial}) to obtain the following key result, which we prove at the
end of the section.
\begin{thm}\label{neweraxial} The monomials $c_iX_1^iX_2^{n-i}$ in the $n$th power of the axial class in
$\tmf^{8n}(RP^\infty\times RP^\infty)$ are equal to those in the $n$th power of
\begin{equation}\label{decomp}
(X_1+X_2)\biggl(u+\sum_{i\ge1}2^{4+i}\a_i\bigl(\bigl(\frac{X_1}{X_2}\bigr)^i+\bigl(\frac{X_2}{X_1}\bigr)^i\bigr)\biggr),
\end{equation}
where $u$ is an odd 2-adic integer and $\a_i$ are 2-adic integers. 
\end{thm}
The factor which accompanies $(X_1+X_2)$ in (\ref{decomp}) is a unit in $\tmf^*(RP^\infty\times RP^\infty)$; we referred to it earlier as $\bu$.
Indeed, its inverse is a series of the same form, obtained by solving a sequence of equations. This justifies the claim in the first paragraph of this
section regarding retrieval of the nonimmersions of \cite[1.1]{BDM}.

We must also observe that restriction to $\tmf^{8\ell}(RP^n\times RP^m)$ of the non-$X_1^iX_2^{\ell-i}$ parts of the basis of $\tmf^{8\ell}(RP^\infty\times RP^\infty)$ cannot cancel the $X_1^iX_2^{\ell-i}$ terms essential for the
nonimmersion. This is proved by noting that these elements such as $L_1X_2^{\ell}$ and $c_4^iL_1X_2^{\ell + i}$
will restrict to a class of the same name in $\tmf^{8\ell}(RP^n\times RP^m)$, and will be 0 there for dimensional
reasons, since $8\ell>n$. 

\begin{pf*}{Proof of Theorem \ref{neweraxial}}
Let $g^*(X)$ denote the axial class as in (\ref{axialsum}). From (\ref{newaxial}) and \ref{c4thm}, the difference
$g^*(X)-(X_1+X_2)$ equals
$$2\sum_{i\ge1}\k_i(X_1^{i+1}+X_2^{i+1})2^{4i}\biggl(\sum_{j\ge0}2^jA_j\biggl(\frac 1{X_1}\bigl(\frac{X_2}{X_1}\bigr)^j+\frac 1{X_2}\bigl(\frac{X_1}{X_2}\bigr)^j\biggr)
\biggr)^i.$$
We let $z=\sqrt{X_1/X_2}$ and $p_j=z^j+z^{-j}$ as in the proof of \ref{c4thm}.

The summand with $i=2t$ becomes
\begin{eqnarray*}&&2\k_i(X_1+X_2)\frac{\sum_sX_1^{2t-s}X_2^s}{X_1^tX_2^t}2^{4i}\biggl(\sum_{j\ge0}2^jA_jp_{2j+1}\biggr)^i\\
&=&2\k_i(X_1+X_2)(p_{2t}+{\cal{L}})2^{4i}\sum_kc_k2^k(p_{2k+i}+{\cal L}).\end{eqnarray*}
Here $k$ is a sum of $j$-values taken from the various
factors in the $i$th power. Also, in $p_{j}+{\cal L}$, $\cal L$ denotes a combination of $p_t$'s with $t<j$.
Noting $(p_{2t}+{\cal L})(p_{2k+i}+{\cal L})=p_{2k+2i}+{\cal L}$, this becomes
\begin{equation}\label{leq}2(X_1+X_2)2^{4i}\sum c'_k2^k(p_{2k+2i}+\cal L).\end{equation}

The argument when $i=2t+1$ is similar but slightly more complicated because $(X_1^{i+1}+X_2^{i+1})$ is not divisible by $(X_1+X_2)$. We obtain
$$2\k_i\frac{X_1^{i+1}+X_2^{i+1}}{(\sqrt{X_1X_2})^{2t+1}}2^{4i}\bigl(\sum_{j\ge0}2^jA_jp_{2j+1}\bigr)^i.$$
For one of the factors of the $i$th power, say the first, we treat $p_{2j+1}$ as $\frac{X_1+X_2}{\sqrt{X_1X_2}}(
p_{2j}+{\cal L})$. The expression then becomes
$$2(X_1+X_2)p_{i+1}2^{4i}\sum c_k2^k(p_{2k+i-1}+{\cal L}),$$
where $k$ is obtained as in the previous case. We again obtain (\ref{leq}).

Thus when $g^*(X)-(X_1+X_2)$ is written as $(X_1+X_2)\sum\b_jp_{2j}$, the coefficient $\b_j$ satisfies
$\nu(\b_j)\ge (j-1)+4+1.$ Here the $(j-1)+4$ comes from the case $i=1$, $k=j-1$ in (\ref{leq}), and the extra $+1$ is the
factor 2 which has been present all along. This yields the claim of (\ref{decomp}).
\end{pf*}

\section{$\tmf$-cohomology of $CP^\infty\times CP^\infty$}\label{CPsec}
In \cite{As1}, \cite{AD}, and \cite{Ann}, it was noted, first by Astey, that the axial class using $BP$ (or $BP\langle2\rangle$)
was $u(X_2-X_1)$, where $u$ is a unit in $BP^*(P^\infty\w P^\infty)$. In this section, we review that argument
and consider the possibility that it might be true when $BP$ is replaced by $\tmf$, which would render
the considerations of the previous section unnecessary. To do this, we calculate $\tmf^*(CP^\infty)$ and
$\tmf^*(CP^\infty\times CP^\infty)$ in positive dimensions. (See Theorems \ref{CPthm} and \ref{CPxCP}.)
Although our conclusion will be that Astey's
$BP$-argument cannot be adapted to $\tmf$, nevertheless these calculations may be of independent interest.

We begin by reviewing Astey's argument. There is a commutative diagram, in which $RP=RP^\infty$ and $CP=CP^\infty$
$$\begin{CD}RP@>d_R>> RP\times RP@>m_R>> RP\\
@VhVV @Vh\times h VV @VhVV\\
CP @>d_C>> CP\times  CP @. CP\\
@. @V 1\times(-1)VV @V 1 VV\\
{} @. CP\times CP @>m_C>> CP\end{CD}$$
The generator $X_R\in BP^2(RP)$ satisfies $X_R=h^*(X)$. We also have that $m_C\circ(1\times(-1))\circ d_C$ is
null-homotopic.
The key fact, which will fail for $\tmf$, is $BP^*(CP\times CP)\approx BP^*[X_1,X_2]$.

The axial class is $m_R^*(X_R)$. It equals $(h\times h)^*(1\times(-1))^*m_C^*(X)$. But 
$$(1\times(-1))^*m_C^*(X)\in\ker(d_C^*).$$
By the above ``key fact," $d_C^*$ is the projection $BP^*[X_1,X_2]\to BP^*[X]$ in which
each $X_i\mapsto X$. The kernel of this projection is  the ideal $(X_2-X_1)$. To see this,
just note that in grading $2n$ a kernel element must be $\sum c_{i}X_1^iX_2^{n-i}$
with $\sum c_{i}=0$, and hence is 
$$\sum_{i<n}c_{i}(X_1^iX_2^{n-i}-X_1^n)=\sum_{i<n}c_iX_1^i(X_2-X_1)\sum X_1^jX_2^{n-i-1-j}.$$

Thus $(1\times(-1))^*m_C^*(X)=(X_2-X_1)u$ for some $u\in BP^*(CP\times CP)$. This $u$ is a unit
by consideration of its reduction to $H^*(-;\Z)$, as in \cite{As1}. Since $h^*(u)$ will then be
a unit in $BP^*(RP\times RP)$ and $h^*(X_i)={X_R}_i$, we obtain the claim about the axial class
being a unit times ${X_R}_2-{X_R}_1$. 

In order to see if there is any chance of adapting this to $\tmf$, we compute $\tmf^*(CP^\infty)$ and
$\tmf^*(CP^\infty\times CP^\infty)$ in positive gradings.
We begin with the relevant Ext calculations. 

Let $\bo=\ext^{*,*}_{A_1}(\zt,\zt)$. Recall that a chart for this is given as in Diagram \ref{bodiag}, extended
with period $(t-s,s)=(8,4)$.
\begin{center}
\begin{minipage}{6.5in}
\begin{diag}\label{bodiag}{$\ext^{*,*}_{A_1}(\zt,\zt)$}
\begin{center}
\begin{picture}(440,100)
\def\mp{\multiput}
\def\elt{\circle*{3}}
\put(98,5){$0$}
\put(138,5){$4$}
\put(178,5){$8$}
\put(95,20){\line(1,0){110}}
\put(100,20){\vector(0,1){80}}
\mp(100,20)(0,10){7}{\elt}
\mp(100,20)(80,40){2}{\line(1,1){20}}
\mp(110,30)(10,10){2}{\elt}
\put(140,50){\vector(0,1){50}}
\mp(140,50)(0,10){4}{\elt}
\put(180,60){\vector(0,1){40}}
\mp(180,60)(0,10){3}{\elt}
\mp(190,70)(10,10){2}{\elt}
\put(230,85){$\cdots$}
\end{picture}
\end{center}
\end{diag}
\end{minipage}
\end{center}

Let $M_{10}$ denote the $A_2$-module $\langle 1,\sq^4,\sq^2\sq^4,\sq^4\sq^2\sq^4\rangle$.
\begin{lem}\label{M10} There is an additive isomorphism
$$\ext_{A_2}^{*,*}(M_{10},\zt)\approx\bo[v_2],$$
where $v_2\in\ext^{1,7}(-)$. \end{lem}

Thus the chart for $\ext^{*,*}_{A_2}(M_{10},\zt)$ consists of a copy of $\bo$ shifted by
$(t-s,s)=(6i,i)$ units for each $i\ge0$.

\begin{pf} There is a short exact sequence of $A_2$-modules
$$0\to \Sigma^7M_{10}\to A_2/\!/A_1\to M_{10}\to 0.$$
This yields a spectral sequence which builds $\ext_{A_2}^{*,*}(M_{10},\zt)$
from $$\bigoplus_{i\ge0}\ext_{A_2}^{*-i,*-7i}(A_2/\!/A_1,\zt).$$
Since $\ext_{A_2}^{*,*}(A_2/\!/A_1,\zt)\approx\bo$, one easily checks that there are
no possible differentials in this spectral sequence.\end{pf}

Let $\bC^m_n=H^*(CP^m_n;\zt)$.
\begin{thm}\label{CPM} There is an additive isomorphism
$$\ext_{A_2}^{*,*}(\bC^\infty_{-\infty},\zt)\approx\bigoplus_{p\in\Z}\Sigma^{8p-2}\bo[v_2].$$
\end{thm}
\ni Of course $\Sigma$ applied to a module or an Ext group just means to increase the $t$-grading by 1.
\begin{pf} There is a filtration of $\bC^\infty_{-\infty}$ with $F_p/F_{p-1}\approx\Sigma^{8p-2}M_{10}$
for $p\in\Z$. We have $\sq^2\iota_{8p-2}=\sq^4\sq^2\sq^4\iota_{8p-10}$. The same argument used in the
last paragraph of the proof of Corollary \ref{batch} works to initiate an inductive
proof of the Ext-isomorphism claimed in the theorem.    
\end{pf}
\begin{cor}\label{DCP} In gradings $(t-s)$ less than $-1$,
$$\ext_{A_2}^{*,*}(\bC^{-2}_{-\infty},\zt)\approx\bigoplus_{p<0}\Sigma^{8p-2}\bo[v_2].$$
\end{cor}
\begin{pf} There is an exact sequence
$$\to\ext_{A_2}^{s-1,t}(\bC_{-1}^\infty,\zt)\to\ext_{A_2}^{s,t}(\bC_{-\infty}^{-2},\zt)\to\ext_{A_2}^{s,t}(\bC_{-\infty}^\infty,\zt)
@>q_*>>\ext_{A_2}^{s,t}(\bC_{-1}^\infty,\zt)\to.$$
The result is immediate from this and \ref{CPM}, since $q_*$ sends the initial tower in $F_0/F_{-1}$ isomorphically
to the initial tower in $\ext_{A_2}(\bC_{-1}^\infty,\zt)$. \end{pf}

The $A$-modules $\bC_1^\infty$ and $\Sigma^2\bC_{-\infty}^{-2}$ are dual. Thus, by \cite[Prop 4]{Seg},
$$\ext_{A_2}^{s,t}(\zt,\bC_1^\infty)\approx\ext_{A_2}^{s,t}(\Sigma^2\bC^{-2}_{-\infty},\zt).$$

There is a ring structure on $\ext_{A_2}^{*,*}(\zt,\bC_1^{\infty})$. We deduce the following result,
which is pictured in Diagram \ref{CPdiag}.
\begin{cor}\label{P1} In $(t-s)$ gradings $\le0$, there is a ring isomorphism
$$\ext^{*,*}_{A_2}(\zt,\bC_1^\infty)\approx \bo[v_2][X],$$
where $X\in\ext^{0,-8}$.\end{cor}
\begin{pf} We apply the duality isomorphism to \ref{DCP}. The multiplicative structure is obtained
from the observation that the powers of the class in $\ext^{0,-8}$ equal the class in $\ext^{0,-8i}$
for each $i>0$.\end{pf}

The Ext groups computed here are the $E_2$-term of the ASS converging to $\tmf^{-*}(CP^\infty)$.
We will consider the differentials in this spectral sequence after performing the Ext calculation
relevant for $\tmf^*(CP^\infty\times CP^\infty)$.

Now we consider $\bC^{-2}_{-\infty}\otimes\bC^{-2}_{-\infty}$. 
Now $x_1$ and $x_2$ denote elements of $H^2(CP;\zt)$. Let $E_2$ denote the exterior subalgebra generated 
by the Milnor primitives of grading 1, 3, and 7. Note that $A_2/\!/E_2$ has a basis with elements of
grading 0, 2, 4, 6, 6, 8, 10, and 12. Finally we note that for any $j\equiv -2$ mod 8 with $j\le -10$, there is a nontrivial $A_2$-morphism $\bC^{-2}_{-\infty}@>\rho>> \Sigma^j\zt$.
\begin{lem}\label{K} Let 
$$K=\ker(\bC^{-2}_{-\infty}\otimes\bC^{-2}_{-\infty}@>\rho>>\bC^{-2}_{-\infty}\otimes\Sigma^{-10}\zt).$$
Let $S$ denote the set of all classes $x_1^{8i-2}x_2^{8j-2}$ with $i\le-1$ and $j\le-2$, together with the classes
$x_1^{8i-2}x_2^{8j+2}$ with $i\le-1$ and $j\le-1$. Then
$K$ is the direct sum of a free $A_2/\!/E_2$-module on $S$ with a single relation $\sq^4\sq^2\sq^4(x_1^{-10}x_2^{-6})=0$.
\end{lem}
\begin{pf} Since the generators of $E_2$ have odd grading, $A_2/\!/E_2$ acts on any element of these
evenly-graded modules. The action of $A_2/\!/E_2$ on $x_1^{-2}x_2^{-2}$ yields the additional elements $x_1^{-2}x_2^0
+x_1^0x_2^{-2}$,
$x_1^{-2}x_2^2+x_1^0x_2^0+x_1^2x_2^{-2}$, $x_1^{-2}x_2^4+x_1^4x_2^{-2}$, $x_1^0x_2^2+x_1^2x_2^0$, $x_1^0x_2^4+x_1^4x_2^0$, 
$x_1^{-2}x_2^8+x_1^2x_2^4+x_1^4x_2^2+x_1^8x_2^{-2}$, and $x_1^0x_2^8+x_1^8x_2^0$.
The action of $A_2/\!/E_2$ on $x_1^{-2}x_2^2$ yields the additional elements $x_1^0x_2^2+x_1^{-2}x_2^4$, $x_1^0x_2^4+
x_1^2x_2^2$, $x_1^2x_2^4+x_1^4x_2^2$,
$x_1^2x_2^4+x_1^{-2}x_2^8$, $x_1^0x_2^8+x_1^4x_2^4$, $x_1^2x_2^8+x_1^8x_2^2$, and $x_1^4x_2^8+x_1^8x_2^4$. 
Each exponent can be decreased by any multiple of 8.

One can easily check that in each grading all classes in $\bC^{-2}_{-\infty}\otimes\bC^{-2}_{-\infty}$
are obtained exactly once from the described elements in $K$ together with $\bC^{-2}_{-\infty}\otimes\Sigma^{-10}\zt$.
There are four cases, for the four even mod 8 values. We illustrate with the case of grading 4 mod 8.
We will just consider the specific value $-28$, but it will be clear that it generalizes to all gradings
$\equiv 4$ mod 8. Letting $X_i$ denote $x_1^ix_2^{-28-i}$, we have:
\begin{enumerate}
\item From generators in $-28$, we obtain just $X_{-10}$ in $K$. The class $X_{-18}$ is in $\bC^{-2}_{-\infty}\otimes\Sigma^{-10}\zt$.
\item From generators in $-32$, we obtain $X_{-8}+X_{-6}$, $X_{-16}+X_{-14}$, and $X_{-24}+X_{-22}$.
\item From generators in $-36$, we obtain $X_{-8}+X_{-4}$ and $X_{-16}+X_{-12}$.
\item From generators in $-40$, we obtain $X_{-4}$, $X_{-12}+X_{-8}$, $X_{-20}+X_{-16}$, and $X_{-24}$.
\end{enumerate}
Note in (4) that $X_0$ and $X_{-28}$ do not appear because each component must be $\le-4$ and the components
sum to $-28$.

One easily checks that the 11 classes listed above, including $X_{-18}$, form a basis for the
space spanned by $X_{-4},\ldots,X_{-24}$, in an orderly fashion that clearly generalizes to
any grading $\equiv4$ mod 8. A similar argument works in the other three congruences.
There are some minor variations in the top few dimensions. \end{pf}

Now we dualize. There is a pairing
$$\ext_{A_2}(\zt,\bC^{\infty}_1)\otimes\ext_{A_2}(\zt,\bC^{\infty}_1)\to\ext_{A_2}(\zt,\bC^{\infty}_1\otimes\bC^{\infty}_1).$$
Let $X_i$  denote the class in grading $-8$ coming from the $i$th factor.
Then we obtain
\begin{thm}\label{PP} The algebra $\ext^{0,*}_{A_2}(\zt,\bC^{\infty}_1\otimes\bC^{\infty}_1)$ in gradings $\le-8$
is isomorphic to $\zt[X_1,X_2]\langle X_1X_2,y_{-12}\rangle$ with $y_{-12}^2=X_1^2X_2+X_1X_2^2$.
The monomials of the form $X_1^iX_2^jy_{-12}$ are acted on freely by $\zt[v_0,v_1,v_2]$. Let
$S_n$ denote the $\zt$-vector space with basis the monomials $X_1^iX_2^{n-i}$, and define a homomorphism
$\eps:S_n\to\zt$ by sending each monomial to 1.
Then $\zt[v_0,v_1,v_2]$ acts freely on $\ker(\eps)$, while $\bo[v_2]$ acts freely on $S_n/\ker(\eps)$.
Thus in dimensions $t-s\le-8$ $\ext^{*,*}_{A_2}(\zt,\bC^{\infty}_1\otimes\bC^{\infty}_1)$ has, for each $i>0$, $i$ copies of $\Sigma^{-8i-4}\zt[v_0,v_1,v_2]$ and $i$ copies of $\Sigma^{-8i-16}\zt[v_0,v_1,v_2]$, and also one copy of
$\Sigma^{-8i-8}\bo[v_2]$.
\end{thm}

Here $\zt[X_1,X_2]\langle X_1X_2,y_{-12}\rangle$ means a free $\zt[X_1,X_2]$-module on basis $\{X_1X_2,y_{-12}\}$
\begin{pf}  The structure as graded abelian group is straightforward from Lemma \ref{K}, Corollary \ref{P1}, and the duality isomorphism
$$\ext^{*,*}_{A_2}(\zt,\bC^{\infty}_1\otimes \bC^{\infty}_1)\approx \ext^{*,*-4}_{A_2}(\bC^{-2}_{-\infty}\otimes\bC^{-2}_{-\infty},\zt).$$
We use that $\ext_{A_2}(A_2/\!/E_2,\zt)\approx\zt[v_0,v_1,v_2]$. The reason that we only assert the
structure in dimension $\le-8$ is due to the $\Sigma^{-10}$ in the cokernel part of Lemma \ref{K}, and
that Theorem \ref{P1} was only valid in dimension $\le0$. In the range under consideration,
the relation on the top class in Lemma \ref{K} does not affect Ext.

The ring structure in filtration 0 comes from $\Hom_{A_2}(\zt,\bC^{\infty}_1\otimes \bC^{\infty}_1)$
being isomorphic to elements of $\bC^{\infty}_1\otimes \bC^{\infty}_1$ annihilated by $\sq^2$ and $\sq^4$,
which has as basis all elements $x_1^{4i}\otimes x_2^{4j}$ and $(x_1^{4i}\otimes x_2^{4j})(x_1^4\otimes x_2^2+x_1^2+x_2^4)$.

Now we show that $\ext_{A_2}^{1,-8n+2}(\zt,\bC^{\infty}_1\otimes \bC^{\infty}_1)=\zt$, and $h_1$ times each monomial in $\ext_{A_2}^{0,-8n}(\zt,\bC^{\infty}_1\otimes \bC^{\infty}_1)$ equals the nonzero element here. An element in $\ext_{A_2}^{1,-8n+2}(\zt,\bC^{\infty}_1\otimes \bC^{\infty}_1)=\zt$ is an equivalence class of morphisms $$\Sigma^2A_2\oplus\Sigma^4A_2@>h>>\bC^{\infty}_1\otimes \bC^{\infty}_1$$
which increase grading by $8n-2$, and yield a trivial composite when preceded by
$$\Sigma^4A_2\oplus\Sigma^8A_2@>\begin{pmatrix}\sq^2&\sq^6\\ 0&\sq^4\end{pmatrix}>>\Sigma^2A_2\oplus\Sigma^4A_2.$$
Morphisms $h$ which can be factored as
\begin{equation}\label{k}\Sigma^2A_2\oplus \Sigma^4A_2@>\sq^2,\sq^4>>A_2@>k>>\bC^{\infty}_1\otimes \bC^{\infty}_1
\end{equation}
are equivalent to 0 in Ext.

We illustrate with the case $n=3$. There are $A_2$-morphisms increasing grading by 22 sending either $\Sigma^2A_2$ or $\Sigma^4A_2$ to any one of the following classes:
\begin{equation}\label{list}x_1^1x_2^{12},\ x_1^2x_2^{10},\ x_1^4x_2^9,\ x_1^4x_2^8,\ x_1^5x_2^8,\ x_1^6x_2^6,\ x_1^8x_2^5,\ x_1^8x_2^4,\ x_1^9x_2^4,\ x_1^{10}x_2^2,\ x_1^{12}x_2^1.\end{equation}
The classes are listed in this order because any two adjacent monomials are equivalent using as
$k$ in (\ref{k}) the morphism sending the generator to the indicated classes in succession:
$$x_1^1x_2^{10},\ x_1^2x_2^9,\ x_1^4x_2^7,\ x_1^3x_2^8,\ x_1^5x_2^6,\ x_1^6x_2^5,\ x_1^8x_2^3,\ x_1^7x_2^4,\ x_1^9x_2^2,\ x_1^{10}x_2^1.$$
For example, $(\sq^2,\sq^4)(x_1^{1}x_2^{10})=(x_1^2x_2^{10},x_1^1x_2^{12})$. Thus all classes in (\ref{list})
are equivalent to one another.

That $h_1$ times any monomial $X_1^iX_2^{n-i}$ equals this nonzero element of $\ext_{A_2}^{1,8n+2}(\zt,\bC^{\infty}_1\otimes \bC^{\infty}_1)$ follows from usual Yoneda product consideration.
If $0\leftarrow \zt\leftarrow C_0\leftarrow C_1\leftarrow$ is the beginning of a minimal $A_2$-resolution, with $C_1=\Sigma^1 A_2\oplus\Sigma^2 A_2\oplus\Sigma^4 A_2$, then $h_1X_1^iX_2^{n-i}$ is represented by the composite
$C_1\to C_0\to \bC^{\infty}_1\otimes \bC^{\infty}_1$ sending $\iota_2\mapsto \iota\mapsto X_1^iX_2^{n-i}$, and this is
equivalent to the element described in the previous paragraph.
\end{pf}

Here is a schematic way of picturing Theorem \ref{PP}. We first list the generators in grading greater than $-32$.
Then for each of the two types of generators, we list the structure arising from them in the first 10 dimensions.
The $\bo[v_2]$-structure in the left half of Diagram \ref{str} arises from one tower in dimensions $-24$ and $-16$, while the $\zt[v_0,v_1,v_2]$-structure in the right half of diagram \ref{str} arises from the other towers in
Diagram \ref{gens}.

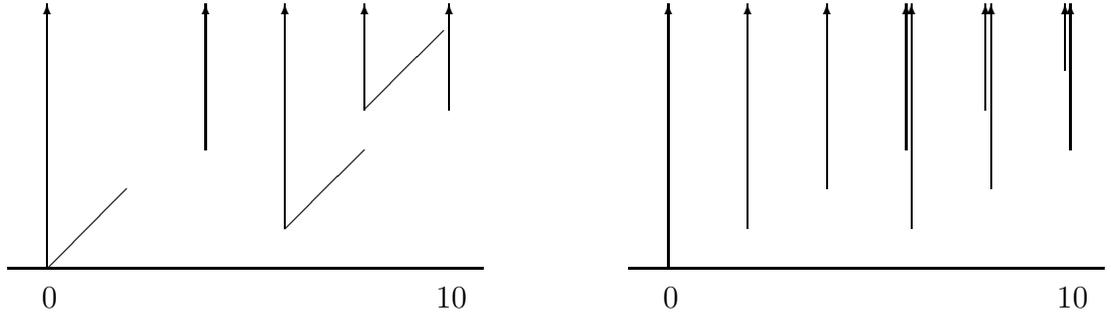
\begin{center}
\begin{minipage}{6.5in}
\begin{diag}\label{gens}{Generators of $\ext_{A_2}(\zt,\bC^\infty_1\otimes \bC^\infty_1)$}
\begin{center}
\begin{picture}(440,130)
\def\mp{\multiput}
\def\elt{\circle*{3}}
\put(100,15){\line(1,0){270}}
\put(105,0){$-28$}
\put(165,0){$-24$}
\put(225,0){$-20$}
\put(285,0){$-16$}
\put(345,0){$-12$}
\mp(115,15)(2,0){3}{\vector(0,1){105}}
\mp(175,15)(2,0){2}{\vector(0,1){105}}
\mp(235,15)(2,0){2}{\vector(0,1){105}}
\mp(295,15)(60,0){2}{\vector(0,1){105}}
\mp(177,15)(118,0){2}{\line(1,1){30}}
\end{picture}
\end{center}
\end{diag}
\end{minipage}
\end{center}

\begin{center}
\begin{minipage}{6.5in}
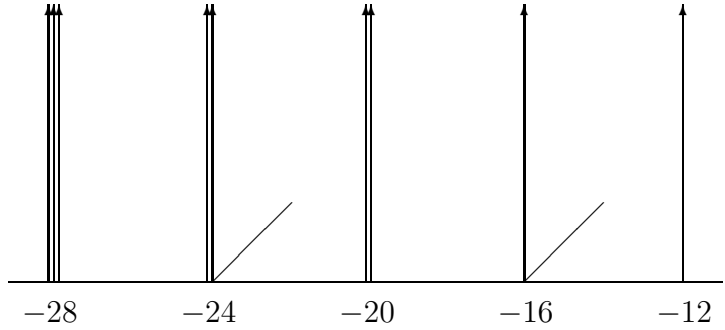
\begin{diag}\label{str}{Structure on two types of generators}
\begin{center}
\begin{picture}(440,130)
\def\mp{\multiput}
\def\elt{\circle*{3}}
\mp(0,15)(235,0){2}{\line(1,0){180}}
\mp(13,0)(235,0){2}{$0$}
\mp(162,0)(235,0){2}{$10$}
\put(15,15){\vector(0,1){100}}
\mp(15,15)(120,60){2}{\line(1,1){30}}
\put(75,60){\vector(0,1){55}}
\put(105,30){\vector(0,1){85}}
\put(105,30){\line(1,1){30}}
\put(135,75){\vector(0,1){40}}
\put(167,75){\vector(0,1){40}}
\put(250,15){\vector(0,1){100}}
\put(280,30){\vector(0,1){85}}
\put(310,45){\vector(0,1){70}}
\put(340,60){\vector(0,1){55}}
\put(370,75){\vector(0,1){40}}
\put(400,90){\vector(0,1){25}}
\put(342,30){\vector(0,1){85}}
\put(372,45){\vector(0,1){70}}
\put(402,60){\vector(0,1){55}}
\end{picture}
\end{center}
\end{diag}
\end{minipage}
\end{center}

\medskip
Now we consider the differentials in the ASS converging to $\tmf^*(CP^\infty)$ and then
for $\tmf^*(CP^\infty\w CP^\infty)$. The gradings are negated when considered as $\tmf$-cohomology
groups.   Corollary \ref{P1}
gives the $E_2$-term converging to $[\Sigma^*CP_1^\infty,\tmf]\approx \tmf^{-*}(CP_1^\infty)$.
We will maintain the homotopy gradings until just before the  end. 
In diagram \ref{CPdiag}, we depict a portion of the $E_2$-term of this ASS in gradings $-16$ to 1.
There are also classes in higher filtration arising from powers of $v_1^4$ and $v_2$ acting on
generators in lower grading.  The elements indicated by $\bullet$'s are involved in differentials,
as explained later.

\begin{center}
\begin{minipage}{6.5in}
\begin{diag}\label{CPdiag}{A portion of $E_2$ for $[\Sigma^*CP^\infty,\tmf]$}
\begin{center}
\begin{picture}(440,240)
\def\mp{\multiput}
\def\elt{\circle*{3}}
\def\dt{\circle*{1}}
\put(0,20){\line(1,0){370}}
\put(10,5){$-16$}
\put(171,5){$-8$}
\put(338,5){$0$}
\put(20,20){\vector(0,1){120}}
\mp(20,20)(160,0){2}{\line(1,1){40}}
\put(180,20){\vector(0,1){160}}
\put(100,80){\vector(0,1){80}}
\put(140,40){\vector(0,1){130}}
\mp(140,40)(160,0){2}{\line(1,1){38}}
\put(300,40){\vector(0,1){170}}
\put(182,100){\vector(0,1){80}}
\mp(182,100)(160,0){2}{\line(1,1){38}}
\put(222,100){\vector(0,1){90}}
\put(258,160){\vector(0,1){40}}
\put(260,80){\vector(0,1){120}}
\put(262,62){\vector(0,1){138}}
\put(262,62){\line(1,1){36}}
\put(302,120){\vector(0,1){90}}
\put(302,120){\line(1,1){38}}
\put(342,100){\vector(0,1){120}}
\put(344,120){\vector(0,1){100}}
\put(346,180){\vector(0,1){40}}
\put(346,180){\line(1,1){36}}
\mp(159,59)(19,19){2}{\elt}
\mp(180,20)(20,20){2}{\elt}
\mp(300,40)(19,19){2}{\elt}
\mp(280,80)(18,18){2}{\elt}
\mp(342,100)(19,19){2}{\elt}
\mp(321,139)(19,19){2}{\elt}
\mp(359,123)(-2,4){8}{\dt}
\end{picture}
\end{center}
\end{diag}
\end{minipage}
\end{center}

We will prove the following key result about differentials in this ASS.
\begin{thm}\label{CPdiffl} The nonzero differentials in the ASS converging to
$[\Sigma^*CP^\infty,\tmf]$, $*<1$, are given by
$$d_2(h_1^\eps v_1^{4i}v_2^jX^{-2k+1})=h_1^{\eps+1}v_1^{4i}v_2^{j+1}X^{-2k}$$
for $\eps=0,1$, $i,j\ge0$, $k\ge1$.\end{thm}
Here $h_1$, $v_1^4$, and $v_2$ have the usual $\ext^{s,t}$ gradings $(s,t)=(1,2)$, $(4,12)$, and
$(1,7)$, respectively. 

Diagram \ref{CPdiag} pictures the situation for $k=1$ and small values of $i$ and $j$.
The elements indicated by $\bullet$'s are involved in the differentials. The resulting 
picture is nicer if the filtrations of all classes built on $X^{-2k+1}$ are increased by 1.
There is a nontrivial extension (multiplication by 2) in dimension $-6$ due to the preceding
differential.
This is equivalent to the way that $bu_*$ is formed from $bo_*$ and $\Sigma^2 bo_*$. We obtain
Diagram \ref{aftershift} from Diagram \ref{CPdiag} after the differentials, extensions, and filtration shift
are taken into account.

\begin{center}
\begin{minipage}{6.5in}
\begin{diag}\label{aftershift}{Diagram \ref{CPdiag} after differentials and filtration shift}
\begin{center}
\begin{picture}(440,240)
\def\mp{\multiput}
\def\elt{\circle*{3}}
\put(0,20){\line(1,0){370}}
\put(10,5){$-16$}
\put(171,5){$-8$}
\put(338,5){$0$}
\put(20,20){\vector(0,1){120}}
\put(20,20){\line(1,1){40}}
\put(180,60){\vector(0,1){120}}
\put(100,80){\vector(0,1){80}}
\put(140,40){\vector(0,1){130}}
\put(300,80){\vector(0,1){130}}
\put(182,100){\vector(0,1){80}}
\put(182,100){\line(1,1){38}}
\put(222,80){\vector(0,1){110}}
\put(258,160){\vector(0,1){40}}
\put(260,100){\vector(0,1){100}}
\put(262,62){\vector(0,1){138}}
\put(302,120){\vector(0,1){90}}
\put(342,140){\vector(0,1){80}}
\put(344,100){\vector(0,1){120}}
\put(346,180){\vector(0,1){40}}
\put(346,180){\line(1,1){36}}
\end{picture}
\end{center}
\end{diag}
\end{minipage}
\end{center}
The regular sequence of towers in the chart beginning in filtration 1 in dimension $-10$ is interpreted
as $v_1^iv_2$, $i\ge0$.

After negating dimensions to switch to cohomology indexing, we obtain the following result, which is immediate from \ref{CPdiffl} after the extensions such as just seen are taken into account.
\begin{thm}\label{CPthm} In positive gradings, there is an isomorphism of graded abelian groups
$$\tmf^*(CP_1^\infty)\approx \Z_{(2)}[Z_{16}](bo^*\oplus v_2\Z_{(2)}[v_1,v_2]).$$
Here $Z_{16}\in \tmf^{16}(CP_1^\infty)$, and $|v_1|=-2$ and $|v_2|=-6$.\end{thm}
Recall that $bo^*=bo_{-*}$ with $bo_*$ as suggested in \ref{bodiag}. Much of the ring structure of $\tmf^*(CP_1^\infty)$ is described in \ref{CPthm}, since
$bo_*$ and $v_2\Z_{(2)}[v_1,v_2]$ are rings, and it is quite clear how to multiply an element in $bo_*$ by one in 
$v_2\Z_{(2)}[v_1,v_2]$. Because of the filtration shift that led to the identification of some of the
classes in $v_2\Z_{(2)}[v_1,v_2]$, we hesitate to make any complete claims about the ring structure.

A complete computation of $\tmf^*(CP^\infty)$ was made in \cite{Bauer}. See there especially Theorem 7.1 and
Diagram 7.1. At first glance, the two descriptions appear quite different, but they seem to be compatible.

\begin{pf*}{Proof of Theorem \ref{CPthm}} We first prove that there is a nontrivial class in $[\Sigma^{-16}CP,\tmf]$ detected in filtration 0.
This is obtained using the virtual bundle $8(H-1)-(H^3-H)$, where $H$ denotes the complex Hopf bundle.
Considered as a real bundle $\theta$, this bundle satisfies $w_2(\theta)$ and $p_1(\theta)=0$. Here we use
from \cite{MS} that $p_1$ generates the infinite cyclic summand in $H^4(BSO;\Z)$ and satisfies $r^*(p_1)=c_1^2-2c_2$
under $BU@>r>> BSO$, and $\rho^*(p_1)=2e_1$ under $BSpin@>\rho>> BSO$, where $H^4(BSpin;\Z)$ is an infinite
cyclic group generated by $e_1$. The total Chern class of $9H-H^3$ is 
$$(1+x)^9(1+3x)^{-1}=1+6x+18x^2+\cdots,$$
and hence 
$$r^*(p_1(\theta))=(c_1(9H-H^3))^2-2c_2(9H-H^3)=(6x)^2-2\cdot18x^2=0.$$
Thus $e_1(\theta)=0$,
hence $CP^\infty@>\theta>>BSpin\to K(\Z,4)$ is trivial, and so $\theta$ lifts to a map $CP^\infty\to BO[8]$.
Hence its Thom spectrum induces a degree-1 map $T(\theta)\to MO[8]$. Since $\psi^3(H)=H^3-H$, by \cite{Sull}
$\theta$ is $J_{(2)}$-equivalent to $8(H-1)$, and hence its Thom spectrum is $T(8(H-1))=\Sigma^{-16}CP_8^\infty$.
Using the Ando-Hopkins-Rezk orientation (\cite{AH}) $MO[8]\to \tmf$, we obtain our desired class as the composite
\begin{equation}\label{AH}\Sigma^{-16}CP_1^\infty@>\text{col}>>\Sigma^{-16}CP_8^\infty@>T(\theta)>>MO[8]\to \tmf.\end{equation}

We will deduce our differentials from the $d_3$-differential $E_3^{4,21}\to E_3^{7,23}$
in the ASS converging to $\pi_*(\tmf)$. This can be seen in \cite[p.537]{Giam} or \cite[Thm 2.2]{DM}.  See
Remark \ref{difflpf} for additional explanation.
It is not difficult to show that, with $M_{10}$ as in \ref{M10}, the morphism
$$\ext_{A_2}^{s,t}(\zt,\zt)\to \ext_{A_2}^{s,t}(M_{10},\zt)$$
induced by the nontrivial $A_2$-map $M_{10}\to\zt$ sends the $\zt$ in $\ext^{7,23}_{A_2}(\zt,\zt)$ which is not
part of the infinite tower to $h_1^2v_1^4v_2$.

We prefer to think about the ASS for $\tmf_*(\Sigma^2CP_{-\infty}^{-2})$, which, as we have noted, is isomorphic
to that of $[\Sigma^*CP_1^\infty,\tmf]$. The $E_2$-term was described in \ref{DCP}. Let $S^{-16}\to \Sigma^2CP_{-\infty}^{-2}\w\tmf$
correspond to the map in (\ref{AH}). Since $E_2(CP_{-\infty}^{-2}\w\tmf)$ in negative dimensions is built from
copies of $\ext_{A_2}(M_{10},\zt)$, we deduce from the previous paragraph that $h_1^2v_1^4v_2g_{-16}$
in the ASS for $\tmf_*(\Sigma^2CP_{-\infty}^{-2})$ must be hit by a $d_2$- or $d_3$-differential, since it is the
image of a class hit by a $d_3$. The only possibility is that it be $d_2$ from $h_1v_1^4g_{-8}$, as indicated
by the dotted line
in Diagram \ref{CPdiag}. Naturality of differentials with respect to $h_1$ and $v_1^4$ implies the differentials
of \ref{CPdiffl} for $\eps=0,1$, all $i$, $j=0$, and $k=1$. Using the diagonal map of $CP_1^\infty$ and the
multiplication of $\tmf$, powers of (\ref{AH}) give similar nontrivial elements in $[\Sigma^{-16k}CP_1^\infty,\tmf]$
for all $k\ge1$,
and by the argument just presented, we establish the differentials of \ref{CPdiffl} for all $k$ (with $j=0$ still).

The only possible differentials on $v_2g_{-16}$ would be some $d_r$ with $r>2$ hitting an element which is
acted on nontrivially by $h_1$. However $h_1v_2g_{-16}$ has become 0 in $E_3$ since it was hit by a $d_2$-differential.
Thus a nonzero differential on $v_2g_{-16}$ would contradict naturality of differentials with respect to
$h_1$-action. Hence there is a map $S^{-10}\to \Sigma^2CP_{-\infty}^{-2}\w\tmf$ hitting $v_2g_{-16}$, and the argument of
the previous paragraph implies that $d_2(h_1v_1^4v_2g_{-8})=h_1^2v_1^4v_2^2g_{-16}$ and then other related
differentials. This now establishes the differentials of \ref{CPdiffl} when $j=1$, and sets in motion an
inductive argument to establish these differentials for all $j\ge1$.

No further differentials in the spectral sequence are possible, by dimensional and $h_1$-naturality considerations. \end{pf*}

\begin{rmk}\label{difflpf}
{\rm The proof of the key $d_3$-differential in the ASS of $\tmf$ from the 17-stem to the 16-stem, which was cited above,
has not had a thorough proof in the literature. Giambalvo's original argument was incorrect and his correction
merely refers to ``a homotopy argument." The current authors cited Giambalvo's result in \cite{DM} without
additional argument. We provide some more detail here regarding this differential.

The relevant portion of the ASS of $\tmf$ appears in Diagram \ref{tmfch}. In \cite{Giam} and \cite{DM},
this was pictured as the ASS of $MO[8]$, but through dimension 18, $$\ext_A^{*,*}(H^*(MO[8]),\zt)\approx \ext^{*,*}_{A_2}
(\zt\oplus\Sigma^{16}\zt,\zt).$$ One way of obtaining the differentials from 15 to 14, as in \cite{Giam}, is to
note that the [8]-cobordism group of 14-dimensional manifolds is $\zt$, and so the top two elements must be killed
by differentials. It is not difficult to compute in Ext the Massey product formula $B=\langle A,h_0,h_1\rangle$,
where $A$ and $B$ are as in Diagram \ref{tmfch}.
This can be seen as $v_1^4$ times a similar formula between classes in dimensions 6 and 8. Since $A$ is 0 in homotopy,
the associated Toda bracket formula says that $B$ must be divisible by $\eta$. But only 0 can be divisible by
$\eta$ in dimension 16 here. Thus $B$ must be killed by a differential, and the depicted way is the only way this
can happen. }
\end{rmk}  

 \begin{center}
 \begin{minipage}{6.5in}
 \begin{diag}\label{tmfch}{Portion of ASS of $\tmf$}   \begin{center}            
  \begin{picture}(440,170)  
  \def\mp{\multiput}        
  \def\elt{\circle*{3}}     
 \mp(125,50)(0,25){3}{\elt}
 \mp(150,25)(0,25){3}{\elt}
 \mp(125,50)(25,-25){2}{\line(0,1){50}}
 \mp(150,25)(0,25){2}{\line(-1,2){25}}
 \put(125,50){\line(1,1){25}}
 \mp(200,50)(0,25){3}{\elt}
 \put(200,50){\line(0,1){50}}
 \mp(225,50)(0,25){2}{\elt}
 \put(225,50){\line(0,1){25}}
 \put(200,50){\line(1,1){25}}
 \mp(200,50)(25,25){2}{\line(-1,3){25}}
 \put(225,50){\line(-1,2){25}}
 \put(175,125){\line(1,1){25}}
 \mp(175,125)(25,25){2}{\elt}
 \put(122,10){$14$}
 \put(172,10){$16$}
 \put(222,10){$18$}
 \put(100,21){$3$}
 \put(100,71){$5$}
 \put(100,146){$8$}
 \put(114,96){$A$}
 \put(164,121){$B$}
 \mp(111,20)(0,145){2}{\line(1,0){129}}
 \mp(111,20)(129,0){2}{\line(0,1){145}}
 \end{picture}
 \end{center}
 \end{diag}
 \end{minipage}
 \end{center}
 
 The differentials in the ASS converging to $\tmf_*(CP^{-2}_{-\infty}\w CP^{-2}_{-\infty})$ 
 are implied by the same considerations that worked for $CP^{-2}_{-\infty}$.
 The $\zt[v_0,v_1,v_2]$-parts in Theorem \ref{PP} cannot support differentials
 by dimensionality and $h_1$-naturality. For the $bo$-like part, we prefer thinking about
 it as $[\Sigma^{*+4}CP_1^\infty\w CP_1^\infty,\tmf]\approx \tmf^{-*-4}(CP_1^\infty\w CP_1^\infty)$, where the product structure is more
 apparent.

 Let $Z_n$ denote the nonzero element of $\ext_{A_2}^{0,-8n}(\zt,\bC_1^\infty\otimes \bC_1^\infty)/\ker(h_1)$.
 By Theorem \ref{PP}, $Z_n$ can be represented by $X_1^iX_2^{n-i}$ for any $1\le i<n$.  If $n$ is even and $n\ge4$, choosing $i$ even, $Z_n$ is an infinite cycle because it is an external product
 of infinite cycles. Hence by the proof of Theorem \ref{CPdiffl},
 $$d_2(h_1^\eps v_1^{4i}v_2^jZ_{2k-1})=h_1^{\eps+1}v_1^{4i}v_2^{j+1}Z_{2k}$$
 for $\eps=0,1$, $i,j\ge 0$, and $k\ge2$.
 
 Finally, $X_1X_2$ is an infinite cycle since there is nothing that it can hit. Also,
 $h_1v_2X_1X_2$ and $h_1^2v_2X_1X_2$ are not hit by differentials since
 $\ext_{A_2}^{0,-8}(\zt,\bC_1^\infty\otimes \bC_1^\infty)=0$ by Theorem \ref{PP}.
 We obtain the following.
 \begin{thm}\label{CPxCP} In grading $\ge10$, there is an isomorphism of graded abelian groups
 $$\tmf^*(CP_1^\infty\w CP_1^\infty)\approx y\Z_{(2)}[v_1,v_2,X_1,X_2]\oplus\bigoplus_{n\ge3}I_n\cdot\Z_{(2)}[v_1,v_2]
 \oplus \Z_{(2)}[Z](bo^*\oplus v_2\Z_{(2)}[v_1,v_2]),$$
 where $|y|=12$, $|X_i|=8$, $|Z|=16$, $|v_1|=-2$, and $|v_2|=-6$. Here
 $I_n=\ker(F_n@>\eps>>\Z)$, where $F_n$ is a free abelian group with basis $\{X_1^iX_2^{n-i}:\,1\le i<n\}$, and 
 $\eps(X_1^iX_2^{n-i})=1$. 
 \end{thm}
 Thus $I_n$ consists of all polynomials of grading $n$ with sum of coefficients equal to 0. We could have
 extended the description in \ref{CPxCP} down to grading 8, but the description would have been slightly more complicated,
 since it would include $h_1v_2Z$ and $h_1^2v_2Z$.
 
 The motivation for this section was to see if perhaps $$\ker(\tmf^*(CP^\infty\times CP^\infty)@>d^*>>\tmf^*(CP^\infty))$$ might
 be something nice like the $I(X_1-X_2)$ which was the case for $BP^*(-)$. In Theorem \ref{CPxCP}, we described
 $\tmf^*(CP^\infty\w CP^\infty)$. To obtain $\tmf^*(CP^\infty\times CP^\infty)$, we add on two copies of
 $\tmf^*(CP^\infty)$, which was described in \ref{CPthm}. Denote by $Z_1$ and $Z_2$ the generators in
 $\tmf^{16}(CP^\infty\times CP^\infty)$. Monomials $Z_1^iZ_2^{n-i}$ should equal $Z^n$ of \ref{CPxCP} plus perhaps
 elements of $I_{2n}$ of \ref{CPxCP}. The class $y$ of \ref{CPxCP} plus perhaps a sum of elements of higher
 filtration is in $\ker(d^*)$ and not in the ideal generated by $(Z_1-Z_2)$. Thus, as expected, $\ker(d^*)$ does 
 not have the nice form that it did for $BP^*(-)$, and so we cannot use this argument to show that the axial class in $\tmf^*(RP^\infty\times
 RP^\infty)$ is $u(X_1-X_2)$. However, we showed something like this by a completely different method in Theorem 
 \ref{neweraxial}. We feel that
 the results obtained in Theorems \ref{CPthm} and \ref{CPxCP} should be
 of independent interest.

\def\line{\rule{.6in}{.6pt}}

\end{document}